%
%

\input amstex  

\input amssym
\input amssym.def

\magnification 1200
\loadmsbm
\parindent 0 cm


\define\nl{\bigskip\item{}}
\define\snl{\smallskip\item{}}
\define\inspr #1{\parindent=20pt\bigskip\bf\item{#1}}
\define\iinspr #1{\parindent=27pt\bigskip\bf\item{#1}}
\define\einspr{\parindent=0cm\bigskip}

\define\ot{\otimes}



\centerline{\bf Weak Multiplier Hopf Algebras I}
\centerline{The main theory}
\bigskip
\centerline{\it  Alfons Van Daele \rm $^{(1)}$ and \it Shuanhong Wang \rm $^{(2)}$}
\bigskip\bigskip
{\bf Abstract} 
\nl 
A {\it weak multiplier Hopf algebra} is a pair $(A,\Delta)$ of a non-degenerate idempotent algebra $A$ and a coproduct $\Delta$ on $A$. The coproduct is a coassociative homomorphism from $A$ to the multiplier algebra $M(A\ot A)$ with some natural extra properties (like the existence of a counit). Further we impose extra but natural {\it conditions on the ranges and the kernels of the canonical maps $T_1$ and $T_2$} defined from $A\ot A$ to $M(A\ot A)$ by
$$T_1(a\ot b)=\Delta(a)(1\ot b)
	\qquad\quad\text{and}\qquad\quad 
	T_2(a\ot b)=(a\ot 1)\Delta(b).$$
The first condition is about the ranges of these maps. It is assumed that there exists an {\it idempotent element} $E\in M(A\ot A)$ such that 
$\Delta(A)(1\ot A)=E(A\ot A)$ and $(A\ot 1)\Delta(A)=(A\ot A)E$.
This element is unique if it exists. Then it is possible to extend the coproduct in a unique way to a homomorphism $\widetilde\Delta:M(A)\to M(A\ot A)$ such that $\widetilde\Delta(1)=E$. In the case of a multiplier Hopf algebra we have $E=1\ot 1$ but this is no longer assumed for weak multiplier Hopf algebras. The second condition determines the behavior of the coproduct on the legs of $E$. We require
$$(\Delta\ot \iota)(E)=(\iota\ot\Delta)(E)=(1\ot E)(E\ot 1)=(E\ot 1)(1\ot E).$$
Finally, the last condition determines the kernels of the canonical maps $T_1$ and $T_2$ in terms of this idempotent $E$ by a very specific relation.
\snl
From these conditions we develop the theory. In particular, {\it we construct a unique antipode} satisfying the expected properties and various other data. Special attention is given to the regular case (that is when the antipode is bijective) and the case of a $^*$-algebra (where regularity is automatic).
\snl
{\it Weak Hopf algebras} are special cases of such {\it weak multiplier Hopf algebras}. Conversely, if the underlying algebra of a (regular) weak multiplier Hopf algebra has an identity, it is a weak Hopf algebra. 
Also any groupoid, {\it finite or not}, yields two weak multiplier Hopf algebras in duality. We will give more and interesting examples of weak multiplier Hopf algebras in [VD-W4] and in [VD-W5]. 
\snl
\nl 
{\it October 2012} (Version 2.1)
\nl
\hrule
\bigskip
\parindent 0.7 cm
\item{($1$)} Department of Mathematics, University of Leuven, Celestijnenlaan 200B,
B-3001 Heverlee, Belgium. {\it E-mail}: Alfons.VanDaele\@wis.kuleuven.be
\item{($2$)} Department of Mathematics, Southeast University, Nanjing 210096, China. {\it E-mail}:  Shuanhwang2002\@yahoo.com {\it or} Shuanhwang\@seu.edu.cn
\parindent 0 cm

\newpage



\bf 0. Introduction \rm
\nl
Let $A$ be an algebra. We will only work with algebras over the field $\Bbb C$ of complex numbers (although most of the results are probably true also for algebras over a more general field). The algebra may or may not have an identity, but we always assume that the product is {\it non-degenerate} (as a bilinear form). Then we can consider the multiplier algebra $M(A)$. It is characterized as the largest algebra containing $A$ as a dense ideal. That $A$ is dense means that $m=0$ if $m\in M(A)$ and if either $ma=0$ for all $a\in A$ or $am=0$ for all $a\in A$. Of course, if already $A$ has an identity, then $M(A)=A$. In a note we wrote on the subject, more information is obtained about this type of 'density' of $A$ in $M(A)$ (see [VD-Ve]).
\snl
An algebra $A$ is called {\it idempotent} if $A^2=A$, that is if any element of $A$ can be written as a sum of products of elements in $A$. All the algebras we consider in this paper are assumed to be idempotent. In fact, we will prove that the underlying algebra of a regular weak multiplier Hopf algebra always has local units (see Proposition 4.9 in Section 4). Then the product is automatically non-degenerate and the algebra is idempotent. 
\snl
The tensor product $A\ot A$ is again a non-degenerate algebra and also here we have the multiplier algebra $M(A\ot A)$. A {\it coproduct} on $A$ is a coassociative homomorphism $\Delta:A\to M(A\ot A)$. It is assumed that the canonical linear maps $T_1$ and $T_2$, defined on $A\ot A$ by 
$$T_1(a\ot b)=\Delta(a)(1\ot b)
	\qquad\quad\text{and}\qquad\quad
		T_2(a\ot b)=(a\ot 1)\Delta(b),$$
have range in $A\ot A$. Then coassociativity of $\Delta$ can be expressed as the commutation rule
$$(\iota\ot T_1)(T_2\ot\iota)=(T_2\ot\iota)(\iota\ot T_1).$$
Here $\iota$ denotes the identity map from $A$ to itself. It is easily seen that this coincides with the original condition of coassociativity if the algebra has a unit. This notion of coassociativity for a coproduct on an algebra without identity was first introduced in [VD1].
\snl
We need some extra properties on the coproduct in order to be sure that {\it its legs cover all of $A$}. This can be guarantied if either there is a counit or if the coproduct is assumed to be full (as we will explain in Section 1).
\snl 
The pair $(A,\Delta)$ is called a {\it multiplier Hopf algebra} if the maps $T_1$ and $T_2$ are bijective from $A\ot A$ to itself. It follows that the coproduct is non-degenerate in the sense that
$$\Delta(A)(A\ot A)=(A\ot A)
	\qquad\quad\text{and}\qquad\quad
		(A\ot A)\Delta(A)=(A\ot A).$$
Then $\Delta$ has a unique extension to a unital homomorphism $\widetilde\Delta:M(A)\to M(A\ot A)$. We often use also $\Delta$ to denote this extension. The homomorphisms $\Delta\ot\iota$ and $\iota\ot\Delta$ are non-degenerate as well and then coassociativity can be written again as
$$(\Delta\ot\iota)\Delta=(\iota\ot\Delta)\Delta$$
where now $\Delta\ot\iota$ and $\iota\ot\Delta$ are the unique extensions of these maps on $A\ot A$ to the multiplier algebra $M(A\ot A)$.
\snl
In the theory of {\it weak multiplier Hopf algebras}, as we study it in this paper, it is no longer assumed that the maps $T_1$ and $T_2$ are bijective. Instead there are very specific conditions on the ranges and the kernels of these maps. This is what we will explain now.
\nl
It is assumed that {\it the ranges of $T_1$ and $T_2$} are given by an idempotent $E\in M(A\ot A)$ in the sense that
$$\Delta(A)(1\ot A)=E(A\ot A)
	\qquad\quad\text{and}\qquad\quad
		(A\ot 1)\Delta(A)=(A\ot A)E.$$
It can easily be shown that this idempotent must be unique. We call it the {\it canonical idempotent} of the weak multiplier Hopf algebra.
\snl
Also under these conditions, it is possible to extend the coproduct to a homomorphism $\widetilde\Delta:M(A)\to M(A\ot A)$ so that $\widetilde\Delta(1)=E$. Again this extension is unique and also here, we will use $\Delta$ to denote this extension. Similarly, the homomorphisms 
$\Delta\ot\iota$ and $\iota\ot\Delta$ have unique extensions to $M(A\ot A)$, now with 
$$(\Delta\ot\iota)(1)=E\ot 1
	\qquad\quad\text{and}\qquad\quad
		(\iota\ot\Delta)(1)=1\ot E,$$
where we use again the same symbols for the extended maps and $1$ both for the identity in $M(A)$ and in $M(A\ot A)$.
\snl
Coassociativity is now also expressed as $(\Delta\ot\iota)\Delta=(\iota\ot\Delta)\Delta$. In fact, this equation is even true on the multiplier algebra $M(A)$. In particular, when applied  on the identity in $M(A)$, this gives the equality
$$(\Delta\ot\iota)(E)=(\iota\ot\Delta)(E).$$
Moreover, it can be shown that the idempotent $(\Delta\ot\iota)(E)$ in $M(A\ot A\ot A)$ is smaller than $E\ot 1$ and $1\ot E$. Recall that for two idempotents $e,f$ we say that $e\leq f$ if $ef=fe=e$.
\snl
Now it seems to be a natural condition to impose that $E\ot 1$ and $1\ot E$ commute and that actually
$$(\Delta\ot\iota)(E)=(E\ot 1)(1\ot E).$$
We refer to our first paper on the subject [VD-W3] where we have motivated the various conditions, in particular the ones above for the idempotent $E$.
\nl
We will now consider the {\it conditions on the kernels of the maps $T_1$ and $T_2$}. We again refer to the preliminary paper [VD-W3] for a discussion and the motivation. As it turns out, in order to get a satisfactory theory of weak multiplier Hopf algebras, in particular with a good antipode, we obtained in [VD-W3] that these kernels are determined already by the idempotent $E$ in a very specific way. This is what we explain now, in a somewhat simplified way. 
\snl
Essentially, there exists a right multiplier $F_1$ of $A\ot A^{\text{op}}$ and a left multiplier $F_2$ of $A^{\text{op}}\ot A$ (where $A^{\text{op}}$ is the algebra $A$ with the opposite product) characterized by the formulas
$$E_{13}(F_1\ot 1)= E_{13}(1\ot E)
	\qquad\quad\text{and}\qquad\quad
		(1\ot F_2)E_{13}=(E\ot 1)E_{13}.\tag"(0.1)"$$
Here, $E_{13}$ stands for the element $E$ as sitting in $M(A\ot A\ot A)$ with its first leg in the first factor and with its second leg in the third factor. More precisely, $E_{13}=(\iota\ot\sigma)(E\ot 1)$ where $\sigma$ is the flip map on $A\ot A$ and $\iota\ot\sigma$ the corresponding map on $A\ot A\ot A$, but extended to  $M(A\ot A\ot A)$. Remark that the elements $F_1$ and $F_2$ are  uniquely determined by these equations.
\snl
It is shown that the maps $F_1$ and $F_2$ are idempotents. We  have that the ranges of the maps $1-F_1$ and $1-F_2$ are in the kernels of  $T_1$ and $T_2$ respectively. More precisely
$$\align T_1((a\ot 1)(1-F_1)(1\ot b))&=0\\
         T_2((a\ot 1)(1-F_2)(1\ot b))&=0
\endalign$$
for all $a,b\in A$. We now use $1$ also for the identity map.
\snl
As a final assumption, we impose the conditions that the kernels of the maps $T_1$ and $T_2$ {\it are precisely equal} to the ranges of these idempotents. That is, we require that
$$\align \text{Ker}(T_1)&=(A\ot 1)(1-F_1)(1\ot A) \\
         \text{Ker}(T_2)&=(A\ot 1)(1-F_2)(1\ot A)
\endalign$$
where we use $\text{Ker}(T_1)$ and $\text{Ker}(T_2)$ for these kernels. Remark that this is indeed an extra assumption. If e.g.\ the maps  $T_1$ and $T_2$ are surjective, then $E=1$ and this will give that also $F_1=F_2=1$. On the other hand, we can not expect that the surjectivity of the maps $T_1$ and $T_2$ automatically gives the injectivity. Examples can be constructed from a semi-group with cancellation (such as the additive semi-group of natural numbers $\Bbb N$).
\snl
As mentioned, the above discussion is {\it a bit simplified}. Only in the regular case, we do have such idempotent elements $F_1$ and $F_2$. In the general case, we work with the maps on $A\ot A$, given by these idempotents. We will explain the problem in Section 2 (and further also in Section 3).
\nl
All of this takes us to {\it the notion of a weak multiplier Hopf algebra} as it is studied in this paper. It is a pair $(A,\Delta)$ of an algebra with a coproduct such that the ranges of the maps $T_1$ and $T_2$ are given by such an idempotent $E$ and so the the kernels of the maps $T_1$ and $T_2$ are given as above by the idempotents $F_1$ and $F_2$, determined by $E$ through the formulas (0.1) above. See Definition 1.14 in Section 1 of this paper for the precise formulation.    
\nl
In our preliminary paper on the subject [VD-W3], we have extensively discussed and motivated this definition. We have explained why it is natural to assume the existence of $E$, giving the ranges of the maps $T_1$ and $T_2$, together with the assumptions that
$(\Delta\ot\iota)(E)=(E\ot 1)(1\ot E)$ and that these last two factors commute. We also have explained why the kernels of the maps $T_1$ and $T_2$ should be given by such idempotents $F_1$ and $F_2$. Finally, we have indicated why they should satisfy the formulas (0.1) in order to get a nice  antipode with the expected properties.
\snl
We advise the reader to have a look at the preliminary paper [VD-W3] as it will help to understand better where all the conditions come from. Moreover, in the present paper, we also have more concise arguments in the proofs of some of our results. Very often, details can be found at various places in the preliminary paper. However it is possible to read and understand this paper without reading the preliminary paper as it does not, strictly speaking, assume any knowledge or results from that paper.
\nl
\it Content of the paper \rm
\nl 
In {\it Section 1} we formulate the main definitions. We start with saying what we really mean by a coproduct $\Delta$ and when we call it full. We give the definition of a counit and discuss the relation with the fullness property. We explain the condition given the idempotent $E$ as above and we show how this can be used to extend the coproduct to the multiplier algebra. We denote this extension still by $\Delta$ and we find that $\Delta(1)=E$. We also use the extensions of $\Delta\ot\iota$ and $\iota\ot\Delta$ to explain the extra conditions on $E$. The rest of the section is mainly devoted to the construction of the idempotents $F_1$ and $F_2$ (or rather the associated maps $G_1$ and $G_2$) as discussed before. All this eventually leads to {\it the definition of a weak multiplier Hopf algebra} as it is studied further in this paper. At the end of this section, we also treat the basic examples associated with a groupoid.
\snl
In {\it Section 2} we start with a weak multiplier Hopf algebra as defined at the end of Section 1 and we develop the theory. This mainly consists of {\it constructing the antipode} $S$ and proving its first properties. We also make some important remarks in this section about coverings. 
\snl
We are able to give {\it another characterization} of weak multiplier Hopf algebras when the antipode is assumed to exist. We use this characterization to show that any  weak Hopf algebra is a weak multiplier Hopf algebra. We also look again at the basic examples coming from a groupoid and use this other characterization to give another argument to show that these basic examples give weak multiplier Hopf algebras.
\snl
The main properties of the antipode however are proven in {\it Section 3}. The antipode is a linear map from $A$ to the multiplier algebra $M(A)$. It is an anti-algebra map in the sense that $S(ab)=S(a)S(b)$ holds in $M(A)$ for all $a,b\in A$. It is also non-degenerate and therefore it has a unique extension to a unital anti-homomorphism from $M(A)$ to itself. It is also an anti-coalgebra map. Indeed, we show that also $\Delta(S(a))=\sigma(S\ot S)\Delta(a)$ for all $a\in A$. This formula makes sense in $M(A\ot A)$ because also $S\ot S$ can be extended. It should be mentioned however that this statement is slightly simplified. See Proposition 3.7 in Section 3 for the precise formulation and comments. We finish in Section 3 with a couple of remarks on some of the difficulties that we encounter in the non-regular case.
\snl
{\it Section 4} is completely devoted to {\it the regular case}. A weak multiplier Hopf algebra is regular if and only if the antipode is a bijective map from $A$ to itself. It turns out that the regular case is much nicer than the general case where some peculiar facts can happen (or at least can not be excluded). This case is much better understood than the general case. We give a collection of formulas that can not be shown in the general case but certainly will also help to better understand the whole theory. Here we can show that the underlying algebra of a regular weak multiplier Hopf algebra automatically has local units. In the general case, we were not able to show this, although it seems very natural and related with the other facts that are not very well understood.
\snl
We have to mention however that in this section, we do not fully prove that a weak multiplier Hopf algebra is regular if and only if the antipode is bijective. We only give the proof here for one direction. The other direction is shown in an appendix ({\it Appendix A}). The reason for this is explained. 
\snl
We consider the involutive case in this section. And we look again at the examples coming from a groupoid. 
\snl
In Section 4, we also show that a regular weak multiplier Hopf algebra $(A,\Delta)$ has to be a weak Hopf algebra if the underlying algebra has an identity. Because a finite-dimensional weak Hopf algebra always has a bijective antipode, this shows that the finite-dimensional weak Hopf algebras are precisely the regular weak multiplier Hopf algebras with a finite-dimensional underlying algebra. We also include a short discussion on the axioms for weak Hopf algebras. 
\snl
In {\it Section 5} we draw some conclusions and we discuss possible further research. In particular, we refer to a number of forthcoming papers on the subject. In {\it Weak Multiplier Hopf Algebras II. The source and target algebras} [VD-W4], we continue with the investigation of the source and target algebras, as well as of the source and target maps. In that paper, we also consider some special examples of weak multiplier Hopf algebras. In {\it Weak Multiplier Hopf algebras III. Integrals and duality} [VD-W5], we study integrals on (regular) weak multiplier Hopf algebras and the notion of duality. 
\nl
\it Conventions, basic notions and references \rm
\nl
As mentioned already, we only work with algebras $A$ over $\Bbb C$. We do not assume that they are unital but we need that the product is non-degenerate. We also assume our algebras to be idempotent. As mentioned, we will show that in the regular case, our algebras will have local units and then non-degeneracy of the algebra, as well as the fact that it is idempotent, are automatic.
\snl
When $A$ is such an algebra, we use $M(A)$ for the multiplier algebra of $A$. We will use $L(A)$ and $R(A)$ to denote the {\it left}, respectively the {\it right multipliers}. We write $ma$ when $m$ is a left multiplier and $a\in A$. By definition, $ma\in A$ and we have $m(ab)=(ma)b$ for all $a,b\in A$. Similarly when $m$ is a right multiplier. If $m\in M(A)$, we have a right and a left multiplier and by definition we have $a(mb)=(am)b$ for all $a,b$.  
\snl
We consider $A\ot A$, the tensor product of $A$ with itself. It is again an idempotent, non-degenerate algebra and we can consider the multiplier algebra $M(A\ot A)$. The same is true for a multiple tensor product. In a short note, we have collected some more results about multiplier algebras and the relation with local units (see [VD-Ve]).
\snl
We use $1$ for the identity in any of these  multiplier algebras. On the other hand, we mostly use $\iota$ for the identity map on $A$ (or other spaces), although sometimes, we also write $1$ for  this map. The identity element in a group is denoted by $e$. 
\snl
If $G$ is a groupoid, we will also use $e$ for units. Units are considered as being elements of the groupoid and we use $s$ and $t$ for the source and target maps from $G$ to the set of units. 
\snl
When $A$ is an algebra, we denote by $A^{\text{op}}$ the algebra obtained from $A$ by reversing the product. When $\Delta$ is a coproduct on $A$, we denote by $\Delta^{\text{cop}}$ the coproduct on $A$ obtained by composing $\Delta$ with the flip map (usually denoted by $\sigma$ and extended to the multiplier algebra).
\snl
The range of a linear map $T$ is denoted by $\text{Ran}(T)$ whereas the kernel of such a map is denoted by $\text{Ker}(T)$. 
\snl
The {\it leg numbering} notation is used for elements and for maps. If e.g.\ $\Delta$ is a coproduct from $A$ to $M(A\ot A)$, we will not really use $\Delta_{12}$ and $\Delta_{23}$ for the maps $\Delta\ot\iota$ and $\iota\ot \Delta$ from $A$ to $M(A\ot A\ot A)$ but we will need $\Delta_{13}$, the map from $A$ to $M(A\ot A\ot A)$ that is e.g.\ defined as $\sigma_{23}\circ(\Delta\ot\iota)$ where $\sigma_{23}$ is the map on $A\ot A\ot A$ that flips the last two factors (and extended to the multiplier algebra). Similarly, if e.g.\ $E$ is an element in $M(A\ot A)$, we can consider elements $E_{12}$, $E_{23}$ and $E_{13}$ in $M(A\ot A\ot A)$. For the first two we have $E_{12}=E\ot 1$ and $E_{23}=1\ot E$, whereas for the third one we have $E_{23}=\sigma_{23}(E_{12})$.
\snl
We will also make use of the {\it Sweedler notation}. For a coproduct $\Delta$, as we define it in Definition 1.1 of this paper, we assume that $\Delta(a)(1\ot b)$ and $(a\ot 1)\Delta(b)$ are in $A\ot A$ for all $a,b\in A$. This allows us to make use of the Sweedler notation for the coproduct. The reader who wants to have a deeper understanding of this, is referred to [VD3] where the use of the Sweedler notation for coproducts that do not map into the tensor product, but rather in its multiplier algebra is explained in detail.
\nl
For the theory of Hopf algebras, we refer to the standard works of Abe [A] and Sweedler [S]. For multiplier Hopf algebras and integrals on multiplier Hopf algebras, we refer to [VD1] and [VD2]. Weak  Hopf algebras have been studied in [B-N-S] and [B-S] and more results are found in [N] and [N-V1]. Various other references on the subject can be found in [V]. In particular, we refer to [N-V2] because we will use notations and conventions from this paper when dealing with weak Hopf algebras.
\snl
For the theory of groupoids, we refer to [Br], [H], [P] and [R]. 
\nl
\bf Acknowledgements \rm
\nl
The first named author (Alfons Van Daele) would like to express his thanks to the following people. He is greatly indebted to Leonid Vainerman for introducing him to the subject of weak Hopf algebras and to his coauthor Shuanhong Wang for motivating him to start the research on weak multiplier Hopf algebras. He is also grateful to Michel Enock and Michel Vallin for information discussions on the topological theory of locally compact quantum groupoids (or measured quantum groupoids).
\snl
The second named author (Shuanhong Wang) would like to thank his coauthor for his help and advices when visiting the Department of Mathematics of the K.U.\ Leuven in Belgium during several periods in 2004, 2006 and 2009. He also is grateful to the analysis research group of the K.U.\ Leuven for the warm hospitality. This work is part of a project, supported by a research fellowship from the K.U.\ Leuven in 2004.
 
 
\newpage



\bf 1. The main definitions \rm
\nl
In this section, {\it we will mainly give definitions}. We start with the notion of a coproduct as we use it in this theory and we finish with the definition of a weak multiplier Hopf algebra.
\snl
Let $A$ be an algebra with a non-degenerate product as explained in the introduction. Also the algebra $A\ot A$ then is non-degenerate and we can consider the multiplier algebras $M(A)$ and $M(A\ot A)$. The identity in these multiplier algebras will always be denoted by $1$. If $A$ is a $^*$-algebra, then also $M(A)$ is a $^*$-algebra. Recall that, in a natural way, we have the inclusions 
$$A\ot A\subseteq M(A)\ot M(A) \subseteq M(A\ot A)$$
and that in general, all these inclusions are strict. Of course, if $A$ already has an identity, the product is non-degenerate and we have $M(A)=A$ as well as $M(A\ot A)=A\ot A$. 
\snl
We also assume that the algebra is idempotent, that is $A^2=A$. Again this is automatic if the algebra has an identity.
\snl
We will see later that most likely, the underlying algebra of a weak multiplier Hopf algebra has local units. If this is the case, the product is automatically non-degenerate and the algebra is idempotent. In fact, we will only be able to show the existence of local units in the regular case. But as there are strong indications that this result is also true for the non-regular case, one might as well assume that all our algebras have local units from the very beginning. We will see that this is the case in all the examples we will consider. Still, we will not require the existence of local units to develop our theory.
\nl
\it The notion of a coproduct \rm
\nl 
In this paper, we work with the following notion of a coproduct.

\inspr{1.1} Definition \rm
A {\it coproduct} (or comultiplication) on the algebra $A$ is a homomorphism $\Delta: A\to M(A\ot A)$ such that 
\snl
i) $\Delta(a)(1\ot b)$ and $(a\ot 1)\Delta(b)$ are in $A\ot A$ for all $a,b\in A$, \newline
ii) $\Delta$ is {\it coassociative} in the sense that
$$(a\ot 1\ot 1)(\Delta \ot \iota)(\Delta(b)(1\ot c))=
      (\iota \ot \Delta)((a\ot 1)\Delta(b))(1\ot 1\ot c)$$
for all $a,b,c\in A$, where $\iota$ is used for the identity map on $A$. 
\snl
In the case of a $^*$-algebra we assume that $\Delta$ is a $^*$-homomorphism.
\snl
The coproduct is called {\it regular} if also  $\Delta(a)(b \ot 1)$ and $(1\ot a)\Delta(b)$ are in $A\ot A$ for all $a,b\in A$.
\hfill$\square$\einspr

If $A$ has an identity, then $\Delta:A\to A\ot A$, condition i) is automatic and ii) is equivalent with the usual condition $(\Delta\ot\iota)\Delta=(\iota\ot\Delta)\Delta$. In general however, condition i) is needed to give a meaning to coassociativity as formulated in condition ii). This notion has been introduced first in [VD1].
\snl
When the coproduct is regular, the opposite coproduct $\Delta^{\text{cop}}$, obtained from $\Delta$ by composing it with the flip, is again a coproduct in the sense of Definition 1.1. Also  a regular coproduct on $A$ is still a coproduct on $A^{\text{op}}$, the algebra $A$ with opposite product. Regularity is  automatic if the algebra is abelian or if the coproduct is coabelian, that is if $\Delta=\Delta^{\text{cop}}$. In the involutive case, regularity of the coproduct is  automatic. Of course, that is the case also when $A$ is unital. 

\inspr{1.2} Notation \rm
We define the {\it canonical maps} $T_1$ and $T_2$ on $A\ot A$ by
$$T_1(a\ot b)=\Delta(a)(1\ot b)
	\qquad\quad\text{and}\qquad\quad
		T_2(a\ot b)=(a\ot 1)\Delta(b).$$
\vskip -0.8 cm
\hfill$\square$\einspr
By assumption, these maps have range in $A\ot A$. Moreover, coassociativity of $\Delta$ can be expressed in terms of the canonical maps as
$$(T_2\ot\iota)(\iota\ot T_1)=(\iota\ot T_1)(T_2\ot\iota).$$
\snl
In the case of a regular coproduct, we define also the maps $T_3$ and $T_4$ on $A\ot A$ by
$$T_3(a\ot b)=(1\ot b)\Delta(a)
	\qquad\quad\text{and}\qquad\quad
		T_4(a\ot b)=\Delta(b)(a\ot 1).$$
The conventions we use here will become clear later. Just remark that the maps $T_3$ and $T_4$ are the maps $T_1$ and $T_2$ when the algebra $A$ is replaced by the algebra $A^{\text{op}}$, obtained from $A$ by reversing the product.
\nl
With this definition of a coproduct, nothing prevents $\Delta$ to be completely trivial and this is of course not very useful. We want the coproduct to cover, in a sense, all of $A$. There are two natural ways to formulate this.
\snl
First we want the existence of a counit.

\inspr{1.3} Definition \rm 
Let $\Delta$ be a coproduct on the algebra $A$ as in Definition 1.1. A linear functional $\varepsilon$ on $A$ is called a {\it counit} if
$$(\varepsilon\ot\iota)(\Delta(a)(1\ot b))=ab
\quad\qquad \text{and} \quad\qquad
(\iota\ot\varepsilon)((c\ot 1)\Delta(a))=ca$$
for all $a,b,c\in A$.
\einspr

In the case of a regular coproduct, we also have 
$$(\varepsilon\ot\iota)((1\ot b)\Delta(a))=ba
\quad\qquad \text{and} \quad\qquad
(\iota\ot\varepsilon)(\Delta(a)(c\ot 1))=ac$$
for all $a,b,c\in A$.
\snl
With only this condition, we can not prove that the counit is unique. We need one more condition on the coproduct. It is another way to require that the legs of $\Delta$ are all of $A$. The notion was first introduced  in [VD-W1], see Definition 1.8 in that paper.

\inspr{1.4} Definition \rm
A coproduct $\Delta$ on an algebra $A$ is called {\it full} if the smallest subspaces $V$ and $W$ of $A$ satisfying 
$$\Delta(A)(1\ot A)\subseteq V\ot A
	\qquad\quad\text{and}\qquad\quad
		(A\ot 1)\Delta(A)\subseteq A\ot W$$
are actually both $A$ itself.
\einspr

When $\Delta$ is full, any element $p\in A$ is a linear combination of elements of the form \newline 
$(\iota\ot\omega)(\Delta(a)(1\ot b))$ with $a,b\in A$ and $\omega\in A'$ and similarly on the other side. See e.g.\ Proposition 1.6 in [VD-W1] (and also Lemma 1.11 in [VD-W3]). 
\snl
If the algebra has an identity and if there is a counit, fullness is automatic. If there is no identity but a counit that is a homomorphism, then fullness is also automatic. However, in the theory of weak multiplier Hopf algebras, we do not expect to have a counit that is a homomorphism. If the coproduct is full and if there is a counit, then this counit is unique. We refer to the discussion about this topic in Section 1 of [VD-W3].
\nl
{\it In what follows}, we assume that $(A,\Delta)$ is a pair of a non-degenerate idempotent algebra $A$ with a full coproduct $\Delta$ such that there is a (unique) counit. We do not assume that the coproduct is regular at this moment.
\nl
\it The ranges of the canonical maps $T_1$ and $T_2$ \rm
\nl
We will now formulate conditions about the ranges $\text{Ran}(T_1)$ and $\text{Ran}(T_2)$ of the canonical maps. We first make the following assumption.

\inspr{1.5} Assumption \rm
We assume that there is an element $E\in M(A\ot A)$ such that $E^2=E$ and
$$E(A\ot A)=T_1(A\ot A) 
	\qquad\quad\text{and}\qquad\quad
		(A\ot A)E=T_2(A\ot A).\tag"(1.1)"$$
\einspr

One can easily show that $E$ is unique. In fact, this will follow from the following result.

\inspr{1.6} Proposition \rm
The element $E$ is the smallest idempotent in $M(A\ot A)$ such that
$$E\Delta(a)=\Delta(a)
	\qquad\quad\text{and}\qquad\quad
		\Delta(a)E=\Delta(a) \tag"(1.2)"$$
for all $a\in A$.

\snl\bf Proof\rm: From the assumptions, it follows that 
$$E(\Delta(a)(1\ot b))=\Delta(a)(1\ot b)$$
for all $a,b\in A$. We can cancel $b$ and obtain that $E\Delta(a)=\Delta(a)$ for all $a$. Similarly $\Delta(a)E=\Delta(a)$ for all $a$.
\snl
Now suppose that $E'$ is another idempotent element in $M(A\ot A)$ that satisfies the formulas (1.2) of the proposition. Then $E'E=E$ because $E(A\ot A)=\Delta(A)(1\ot A)$ and $E'\Delta(a)=\Delta(a)$ for all $a$. Similarly $EE'=E$. This means that $E\leq E'$.
\hfill$\square$\einspr

Observe that indeed, this will imply that $E$ is the unique idempotent that satisfies the equations (1.1) in the assumption. We will sometimes call $E$ the {\it canonical idempotent} of the weak multiplier Hopf algebra $(A,\Delta)$.

\inspr{1.7} Remark \rm
In the $^*$-algebra case, where it is assumed that $A$ is a $^*$-algebra and that $\Delta$ is a $^*$-homomorphism, it follows from the above uniqueness that we must have $E^*=E$. In the more general regular case, as treated in detail in Section 3, we will also have that $E$ is the idempotent obtained from $\Delta^{\text{cop}}$ instead of $\Delta$. This will also follow from the proposition above.
\hfill$\square$\einspr

In the case that $E=1\ot 1$, it follows that $\Delta$ is non-degenerate. Then we know that we can extend $\Delta$ to the multiplier algebra and that for the extension, we get $\Delta(1)=1\ot 1$ (see e.g.\ [VD1]). We now generalize this result (see also the appendix in [VD-W3]).

\inspr{1.8} Proposition \rm
There is a unique homomorphism $\widetilde\Delta:M(A)\to M(A\ot A)$ that extends $\Delta$ and so that $\widetilde\Delta(1)=E$.

\snl\bf Proof\rm:
Suppose that $\widetilde\Delta$ is such an extension of $\Delta$ to $M(A)$. Then for all $m\in M(A)$ and $x\in A\ot A$, we will have
$$\widetilde\Delta(m)x
	=\widetilde\Delta(m)\widetilde\Delta(1)x
	=\widetilde\Delta(m)Ex
$$
and this will be equal to $\sum_i \Delta(ma_i)(1\ot b_i)$
where $Ex$ is written as $\sum_i \Delta(a_i)(1\ot b_i)$ with finitely many elements $a_i,b_i\in A$. Such elements exist by assumption. A similar formula is obtained for $y\widetilde\Delta(m)$ if $y\in A\ot A$.
\snl
Now, it is clear that not only the uniqueness of the extension will follow, but also that these formulas can be used to define  $\widetilde\Delta(m)\in M(A\ot A)$ for any $m\in M(A)$.
\snl
Finally, it is straightforward to verify that $\widetilde\Delta$, defined in this way, is a homomorphism, that it extends $\Delta$ and that $\widetilde\Delta(1)=E$.
\hfill$\square$\einspr

As usual, we denote the extension still by $\Delta$.
\snl
In a completely similar way, we can extend the homomorphisms $\Delta\ot\iota$ and $\iota\ot\Delta$ defined from $A\ot A$ to $M(A\ot A\ot A)$  to homomorphisms from  $M(A\ot A)$ to $M(A\ot A\ot A)$ in such a way that 
$$(\Delta\ot\iota)(1)=E\ot 1
	\qquad\quad\text{and}\qquad\quad
		(\iota\ot\Delta)(1)=1\ot E.$$
Here, we use $1$ both for the identity in $M(A)$ and in $M(A\ot A)$ and again we are also using the same symbols for the extensions. Coassociativity is now written in the usual way. In fact, we get 
$$(\Delta\ot\iota)\Delta(m)=(\iota\ot\Delta)\Delta(m)\tag"(1.3)"$$ 
for any element $m\in M(A)$. See e.g. Proposition A.7 in [VD-W3] for details.
\snl
In particular, this is true for the identity in $M(A)$. This gives the following result.

\inspr{1.9} Proposition \rm
We have $(\Delta\ot\iota)(E)=(\iota\ot\Delta)(E)$ and this idempotent in $M(A\ot A\ot A)$ is smaller than both $E\ot 1$ and $1\ot E$.

\snl\bf Proof\rm:
If we apply (1.3) with $m=1$ we find the first formula. And if we use that $E=(1\ot 1)E=E(1\ot 1)$ (with $1\in M(A)$) and apply $\Delta\ot \iota$ we get
$$(\Delta\ot \iota)(E)=(E\ot 1)((\Delta\ot \iota)(E))=((\Delta\ot \iota)(E))(E\ot 1).$$
This precisely means that $(\Delta\ot \iota)(E)\leq (E\ot 1)$. Similarly we get that \newline
$(\iota\ot\Delta)(E)\leq (1\ot E)$.
\hfill$\square$\einspr 

We now need to add the following assumption.

\iinspr{1.10} Assumption \rm
We assume that the idempotents $E\ot 1$ and $1\ot E$ commute and that 
$$(\Delta\ot \iota)(E)=(E\ot 1)(1\ot E).\tag"(1.4)"$$
\einspr

This assumption is quite natural because of the result in Proposition 1.9 above. One may wonder even if it is automatically fulfilled. In the involutive case, the equality (1.4) implies in fact that the idempotents $E\ot 1$ and $1\ot E$ commute because by the equation, the product of them is self-adjoint. For a detailed discussion about this condition, we refer to Section 3 in the preliminary paper [VD-W3] where this is explained and where various reasons are given to make this extra assumption plausible. 
\nl
So far about the conditions on the {\it ranges} of the canonical maps. Now, we turn our attention to the {\it kernels}. 
\nl
\it The kernels of the canonical maps $T_1$ and $T_2$ \rm
\nl
These kernels will also be given by idempotent maps from $A\ot A$ to itself. In the following proposition, we define these maps. Later we show that they have the right properties. In the proof of this proposition, we are using the Sweedler notation for the first time. The reader is advised to verify that the necessary coverings are present to make formulas precise. See a remark in the introduction (and if necessary [VD3]).

\iinspr{1.11} Proposition \rm
There exists  linear maps  $G_1$ and $G_2$ from  $A\ot A$ to itself characterized by the equalities
$$\align (G_1\ot \iota)(\Delta_{13}(a)(1\ot b\ot c))&= \Delta_{13}(a)(1\ot E)(1\ot b\ot c)\\
	(\iota\ot G_2)((a\ot b\ot 1)\Delta_{13}(c))&=(a\ot b\ot 1)(E\ot 1)\Delta_{13}(c)
\endalign$$
for all $a,b,c\in A$. We use the leg-numbering notation as explained in the introduction.

\snl\bf Proof\rm:
We will give the proof for $G_1$. The proof for $G_2$ is completely similar.
\snl
For any $a\in A$ we first define a right multiplier $G(a)$ of $A\ot A$ by the formula
$$(b\ot c)G(a)=\sum_{(a)} ba_{(1)}\ot (\iota\ot\varepsilon)((c\ot a_{(2)})E)$$
where $b,c\in A$. Here $\varepsilon$ is the counit. It is clear that the right hand side of this equation is well-defined in $A\ot A$ and that the equation really defines a right multiplier of $A\ot A$.
\snl
Now take $a$, $b$, $c$ and $d$ in $A$. Then we have
$$\sum_{(a)} (b\ot c)G(a_{(1)}) \ot  a_{(2)}d 
	= \sum_{(a)} (ba_{(1)}\ot (\iota\ot\varepsilon)((c\ot a_{(2)})E)\ot a_{(3)}d. \tag"(1.5)"$$
For any $p\in A$ we have 
$$\align (c\ot \Delta(p))(E\ot 1)(1\ot 1\ot d)
	&=((\iota\ot\Delta)(c\ot p))(1\ot E)(E\ot 1)(1\ot 1\ot d) \\
	&=((\iota\ot\Delta)(c\ot p))(\iota\ot\Delta)(E))(1\ot 1\ot d) \\
	&=((\iota\ot\Delta)(c\ot p)E))(1\ot 1\ot d)
\endalign$$
and if we now apply $\varepsilon$ on the second factor, we get $(c\ot p)E(1\ot d)$. We use this result in the previous formula (1.5), with $p$ replaced by $a_{(2)}$. This gives
$$\sum_{(a)} (b\ot c)G(a_{(1)}) \ot  a_{(2)}d 
	=\sum_{(a)}(ba_{(1)}\ot c\ot a_{(2)})(1\ot E)(1\ot 1\ot d).$$
Now, we multiply with $1\ot q\ot 1$ from the right and we cancel $b$ and $c$. This gives the equality
$$\sum_{(a)} G(a_{(1)})(1\ot q) \ot  a_{(2)}d 
	=\Delta_{13}(a)(1\ot E)(1\ot q\ot d). \tag"(1.6)"$$
Observe that the right hand side is an element in $A\ot A\ot A$ and so is also the left hand side. 
\snl
Finally we use that $\Delta$ is full. This means that any element $p\in A$ is a linear combination of elements of the form $(\iota\ot\omega)(\Delta(a)(1\ot d))$ with $a,d\in A$ and $\omega\in A'$. It follows from this that $G(p)(1\ot q) \in A\ot A$ for all $p,q\in A$. 
\snl
So, we can define the map $G_1:A\ot A\to A\ot A$ by $G_1(p\ot q)=G(p)(1\ot q)$. Moreover, from the formula (1.6) above, we see that 
$$\sum_{(a)} G_1(a_{(1)}\ot q) \ot  a_{(2)}d 
	=\Delta_{13}(a)(1\ot E)(1\ot q\ot d)$$
for all $a,q,d\in A$. This completes the proof.
\hfill$\square$\einspr

It is not hard to prove that  
$$\align G_1(a\ot bb')&=(G_1(a\ot b))(1\ot b')\\
	G_2(a'a\ot b)&=(a'\ot 1)(G_2(a\ot b))
\endalign$$
for all $a,a',b,b'\in A$. We also can show easily that 
$$\align (\Delta\ot\iota)G_1&=(\iota \ot G_1)(\Delta\ot\iota)\\
	(\iota\ot\Delta)G_2&=(G_2\ot\iota)(\iota\ot\Delta).
\endalign$$
The last formulas get meaningful if we multiply with an element $c\in A$ from the left in the first factor (in the first case) and with an element $d\in A$ from the right in the second factor (in the second case). Then in fact, these formulas are read as 
$$\align (T_2\ot\iota)(\iota\ot G_1)&=(\iota \ot G_1)(T_2\ot\iota)\\
	(\iota\ot T_1)(G_2\ot\iota)&=(G_2\ot\iota)(\iota\ot T_1).
\endalign$$
These properties will be crucial in the next section where we use the generalized inverses of $T_1$ and $T_2$, determined by the idempotents that project on the ranges and the kernels, to produce the antipode. For this reason, we will give an explicit proof of these results in the next section (see Proposition 2.2).

\iinspr{1.12} Remark \rm
i) From the considerations in [VD-W3], we expect that there is
a right multiplier $F_1$ of $A\ot A^{\text{op}}$ and a left multiplier $F_2$ of $A^{\text{op}}\ot A$ (where $A^{\text{op}}$ is the algebra $A$ with the opposite product) so that
$$G_1(a\ot b)=(a\ot 1)F_1(1\ot b)
\quad\quad\text{and}\quad\quad
G_2(a\ot b)=(a\ot 1)F_2(1\ot b) \tag"(1.7)"$$
for all $a,b\in A$ and that these multipliers are idempotents. We can show that this is the case if the weak multiplier Hopf algebras is regular (as already briefly considered at the end of this section and properly defined and studied in Section 4). Unfortunately, we are not able to show this in general.  It is possible if we  have that the coproduct is regular (as defined in Definition 1.1) and if also 
$$(1\ot A)\Delta(A)=(A\ot A)E
\qquad\quad\text{and}\quad\qquad
\Delta(A)(A\ot 1)=E(A\ot A).$$
Already under these extra conditions, there exist multipliers $F_1$ and $F_2$ as above, giving the maps $G_1$ and $G_2$ with the formulas (1.7). This will be explained in Section 4 (see a remark following Proposition 4.5).
\snl
ii) In the case of a $^*$-algebra, when $\Delta$ is assumed to be a $^*$-homomorphism, we know that the coproduct is automatically regular and as we have seen that $E=E^*$, the above extra conditions  are automatically fulfilled. We will see in Section 4 that this is also the case for regular weak multiplier Hopf algebras (see Proposition 4.2). In general however, we can not show this. Apparently, this is one of the 'subtleties' of the non-regular case that are not very well understood.
\hfill$\square$\einspr

It would make sense to assume the existence of these multipliers $F_1$ and $F_2$ further in our treatment, but we will not do this as it turns out that we go a long way without this extra assumption. 
\snl
Indeed, in the general case and without this extra assumption, it is possible to prove the following expected property of the maps $G_1$ and $G_2$. In the formulation, we use $1$ (and not $\iota$ as usual) for the identity map on $A\ot A$.

\iinspr{1.13} Proposition \rm
The maps $G_1$ and $G_2$ are idempotent. The ranges of $1-G_1$ and $1-G_2$ are in the kernels of $T_1$ and $T_2$ respectively.

\snl\bf Proof\rm: Again, we only give the proof for $G_1$.
For all $a,b,c\in A$ we have
$$\align (G_1\ot\iota)(G_1\ot\iota)(\Delta_{13}(a)(1\ot b\ot c))
	&=(G_1\ot\iota)(\Delta_{13}(a)(1\ot E)(1\ot b\ot c)) \\
	&=\Delta_{13}(a)(1\ot E)(1\ot E)(1\ot b\ot c) \\
	&=\Delta_{13}(a)(1\ot E)(1\ot b\ot c) \\
	&=(G_1\ot\iota)(\Delta_{13}(a)(1\ot b\ot c)).
\endalign$$
This proves that the map is idempotent. To prove the second statement, take again $a,b,c\in A.$ Then
$$\align (T_1\ot\iota)(G_1\ot\iota)(\Delta_{13}(a)(1\ot b\ot c))
	&=(T_1\ot\iota)(\Delta_{13}(a)(1\ot E)(1\ot b\ot c)) \\
        &=(\Delta\ot\iota)\Delta(a)(1\ot E)(1\ot b\ot c)\\
        &=(\iota\ot\Delta)\Delta(a)(1\ot E)(1\ot b\ot c)\\
        &=(\iota\ot\Delta)\Delta(a)(1\ot b\ot c)\\
        &=(\Delta\ot\iota)\Delta(a)(1\ot b\ot c)\\
 	&=(T_1\ot\iota)(\Delta_{13}(a)(1\ot b\ot c)). 
\endalign$$ 
\vskip -0.7 cm \hfill$\square$\einspr

In the case where $G_1$ and $G_2$ are given by a right multiplier $F_1$ of $A\ot A^{\text{op}}$ and a left multiplier $F_2$ of $A^{\text{op}}\ot A$ as in the formulas (1.7), we will of course have that these multipliers are also idempotents. 


\nl
\it The main definition and some basic examples \rm
\nl
This leads us to the final assumptions where we fix the kernels of the canonical maps $T_1$ and $T_2$ by requiring that they actually coincide with the ranges of the maps $1-G_1$ and $1-G_2$.
\snl
And so, we are ready to formulate the main definition.

\iinspr{1.14} Definition \rm
A {\it weak multiplier Hopf algebra} is a pair $(A,\Delta)$ of a non-degenerate idempotent algebra $A$ with a full coproduct and a counit satisfying the following conditions.
\snl
i) There exists an idempotent $E\in M(A\ot A)$ giving the ranges of the canonical maps: 
$$E(A\ot A)=T_1(A\ot A)
	\qquad\quad\text{and}\qquad\quad
		(A\ot A)E=T_2(A\ot A).$$
ii) The element $E$ satisfies 
$$(\iota\ot \Delta)(E)=(E\ot 1)(1\ot E)=(1\ot E)(E\ot 1).$$\newline
iii) The kernels of the canonical maps are of the form 
$$\align \text{Ker}(T_1)&=(1-G_1)(A\ot A) \\
         \text{Ker}(T_2)&=(1-G_2)(A\ot A),
\endalign$$
where $G_1$ and  $G_2$ are the linear maps from $A\ot A$ to itself, given as in Proposition 1.11.
\hfill$\square$\einspr

We will call $E$ the {\it canonical idempotent} and we will sometimes refer to it as $\Delta(1)$.
\snl
A weak multiplier Hopf algebra will be called {\it regular} if the coproduct is regular and if also $(A,\Delta^{\text{cop}})$ is a weak multiplier Hopf algebra. This is the same as requiring that also $(A^{\text{op}},\Delta)$ is a weak multiplier Hopf algebra. We refer to Section 4 for this case.
\snl
If $A$ is a $^*$-algebra and if $\Delta$ is a $^*$-homomorphism, regularity will be automatic (see Proposition 4.11 in Section 4). 
We call it a {\it weak multiplier Hopf $^*$-algebra}.
\nl
We will now discuss the two basic examples associated with a groupoid. Later, in Section 2 (and also in Section 4), we will consider the case of a (finite-dimensional) weak Hopf algebra and see how this fits into our theory. We will give more examples in [VD-W4] and [VD-W5].

\iinspr{1.15} Example \rm
Let $G$ be any groupoid. Consider the $^*$-algebra $K(G)$ of complex functions with finite support (and pointwise operations). The involution is defined by $f^*(p)=\overline{f(p)}$ when $f\in K(G)$ and $p\in G$. We will denote this $^*$-algebra by $A$. Observe that it only has an identity when $G$ is finite but that the product in $A$ is always non-degenerate and that the algebra is idempotent. The product in $G$ yields a coproduct $\Delta$ on $A$ by the formula
$$
\Delta(f)(p,q)=
\cases
f(pq) & \text{if $pq$ is defined},\\
0 & \text{otherwise}.
\endcases
$$
It is a $^*$-homomorphism from $A$ to the multiplier algebra $M(A\ot A)$ of the tensor product $A\ot A$ of $A$ with itself. Recall that in this case we have a natural identification of $A\ot A$ with $K(G\times G)$ and of $M(A\ot A)$ with $C(G\times G)$, the algebra of all complex functions on $G\times G$. 
\snl 
Take $f,g\in A$ and consider the function $\Delta(f)(1\ot g)$. It maps the pair $(p,q)$ to $f(pq)g(q)$ if $pq$ is defined and to $0$ otherwise. The presence of $g$ forces $q$ to lie in a finite set (for the result to be non-zero). Also $pq$ must lie in a finite set and because $p=(pq)q^{-1}$ when $pq$ is defined, the result will be $0$ except when also $p$ lies in a finite set. Therefore $\Delta(f)(1\ot g)\in K(G\times G)$. Similarly for $(f\ot 1)\Delta(g)$ so that condition i) in Definition 1.1 is satisfied. 
\snl
The coassociativity of $\Delta$, as formulated in ii) of the definition, is a straightforward consequence of the associativity of the product in $G$. Regularity of the coproduct here follows automatically because the algebra is abelian.
\snl
The counit is given by the formula $\varepsilon(f)=\sum f(e)$ where the sum is taken over all the units of $G$. To show that $\Delta$ is full, take any element $p\in G$ and let $e$ be the target of $p$. Take $f=\delta_p$, the function that is $1$ in $p$ and $0$ everywhere else. Then $\Delta(f)(e,\,\cdot\,)=f$ and this shows that the right leg of $\Delta$ is all of $A$. One could also argue that $(\delta_e\ot 1)\Delta(\delta_p)=\delta_e\ot\delta_p$ for all $p$. Similarly, by taking the source of an element, we get that the left leg of $\Delta$ is all of $A$.
\snl
Define $E$ as the function on $G\times G$ that is one on pairs $(p,q)$ for which $pq$ is defined and $0$ on other pairs. Then it is an idempotent element in $M(A\ot A)$ satisfying the condition i) and ii) in Definition 1.14. Indeed, $E$ is the obvious candidate for $\Delta(1)$.
\snl
The maps $G_1$ and $G_2$ are given by idempotents $F_1$ and $F_2$ in $M(A\ot A)$ because the algebra is abelian. These idempotents are functions on $G\times G$. The element $F_1$ is given by the function that is $1$ on pairs $(p,q)$ for which $s(p)=s(q)$ while $F_2$ is given by the function that is $1$ on pairs $(p,q)$ for which $t(p)=t(q)$. Recall that $s(p)$ is used for the source of $p$ and that $t(p)$ denotes the target of $p$. It can be verified that the elements $F_1$ and $F_2$ satisfy the right conditions. Let us look e.g.\ at the formulas, defining the maps $G_1$ and $G_2$ in Proposition 1.11. For the first formula in Proposition 1.11, we get in this case
$(E_{13}(F_1\ot 1))(p,q,v)=1$ if and only if $s(p)=t(v)$ and $s(p)=s(q)$. On the other hand we get $(E_{13}(1\ot E))(p,q,v)=1$ if and only if $s(p)=t(v)$ and $s(q)=t(v)$. These conditions are the same. This proves that $G_1$ is given by this function $F_1$. Similarly, we will find that $G_2$ is given by the function $F_2$. A straightforward argument shows that they provide the kernels of the canonical maps as in the definition of a weak multiplier Hopf algebra. 
\snl 
Needless to add that this gives a regular weak multiplier Hopf algebra (as the algebra is abelian).
\hfill$\square$\einspr  

At various places in the first paper [VD-W3], some more detailed arguments for the above example can be found. In particular, see Proposition 4.3 in [VD-W3] where it is shown that the maps $G_1$ and $G_2$, given by the elements $F_1$ and $F_2$ as above, give the kernels of the canonical maps $T_1$ and $T_2$ respectively. All of this is relatively straightforward. Moreover, at the end of the next section, we will come back to this example and show, in another way, that it fits into our theory. 
\snl
The next example is dual to the previous one.

\iinspr{1.16} Example \rm Again $G$ is a groupoid. Now we take for $A$ the algebra $\Bbb C G$ of complex functions with finite support on $G$ with the {\it convolution product}. If we use $p\mapsto \lambda_p$ for the canonical imbedding of $G$ in $\Bbb C G$, we have $\lambda_p\lambda_q=\lambda_{pq}$ for all $p,q$ so that $pq$ is defined. Otherwise this product is $0$. The algebra $A$ becomes a $^*$-algebra if we let $\lambda_p^*=\lambda_{p^{-1}}$ for all $p$. The algebra has an identity if and only if the set of units in $G$ is finite. In that case, we have $1=\sum \lambda_e$ where the sum is taken over all the units $e$. In the general case, we still have a non-degenerate product and the same formula gives the identity in the multiplier algebra $M(A)$. Remark that the elements $(\lambda_e)$ are mutually orthogonal idempotents in $A$.
\snl
A coproduct on $A$ is given by $\Delta(\lambda_p)=\lambda_p\ot \lambda_p$ for all $p$. It is a $^*$-homomorphism. It maps $A$ to $A\ot A$, also in the case where there is no identity. A counit is given by the usual formula $\varepsilon(\lambda_p)=1$ for all $p$. The coproduct is automatically full. 
\snl
Let $E$ be the element in $M(A\ot A)$ defined as $\sum \lambda_e\ot\lambda_e$ where the sum is again taken over all units. One verifies that $E$ is the smallest idempotent so that
$$E(\lambda_p\ot\lambda_p)=\lambda_p\ot\lambda_p
\qquad\quad\text{and}\quad\qquad
(\lambda_p\ot\lambda_p)E=\lambda_p\ot\lambda_p$$
for all $p$. Clearly also 
$$(\Delta\ot\iota)(E)=(\iota\ot\Delta)(E)=\sum_e \lambda_e\ot \lambda_e \ot \lambda_e.$$
We also have 
$$E\ot 1 = \sum_{e,f} \lambda_e\ot \lambda_e \ot \lambda_f 
\qquad\quad\text{and}\quad\qquad
1\ot E = \sum_{e,f} \lambda_f\ot \lambda_e \ot \lambda_e,
$$ 
where we take the sum over pairs of units, 
and we see that these two elements commute and that the product is precisely $\sum_e \lambda_e\ot \lambda_e \ot \lambda_e$. We use that the elements $(\lambda_e)$ are mutually orthogonal idempotents.
\snl
Further, one verifies that the idempotents $F_1$ and $F_2$, giving the maps $G_1$ and $G_2$ as in Proposition 1.11, are again $\sum_e \lambda_e\ot \lambda_e$. A straightforward argument also shows that they provide the kernels of the canonical maps as in the definition of a weak multiplier Hopf algebra.
\snl
Also here we get a regular weak multiplier Hopf algebra, now because it is coabelian.

\hfill$\square$\einspr

Again, at various places in [VD-W3] (see e.g. Proposition 4.4 in [VD-W3]), more details can be found, but also for this case, the arguments are all straightforward.
\snl
Moreover, at the end of the next section, we will also consider this example again, using the extra results we have found then. As mentioned already, we then will also treat the weak Hopf algebras and show how they fit in our theory.

\nl\nl




\bf 2. The antipode: construction and first properties\rm
\nl
In this section, we start with a weak multiplier Hopf algebra as defined in the previous section (see Definition 1.14). We will first construct '{\it generalized inverses}' $R_1$ and $R_2$ of $T_1$ and $T_2$ respectively. They will yield antipodes $S_1$ and $S_2$ in very much the same way as in the case of ordinary multiplier Hopf algebras (where the maps $T_1$ and $T_2$ are assumed to be bijections). We will see how the specific choices of the projection maps $G_1$ and $G_2$ on the kernels of the canonical maps give the equality $S_1=S_2$ of these antipodal maps, thus defining the antipode $S$ of the weak multiplier Hopf algebra. In the next section, it will be shown that this antipode is an anti-algebra as well as an anti-coalgebra map.
\snl
In this section, we will also give an {\it alternative definition} of a weak multiplier Hopf algebra when an antipode is given. This will be used to show that any (finite-dimensional) weak Hopf algebra is a weak multiplier Hopf algebra. Also the examples coming from a groupoid will be reconsidered from this other point of view.
\nl
So let $(A,\Delta)$ be a weak multiplier Hopf algebra. It is a pair of a non-degenerate, idempotent algebra $A$ with a full coproduct $\Delta$ and a counit $\varepsilon$. There is an idempotent $E$ in $M(A\ot A)$ with the property that the ranges of the canonical maps $T_1$ and $T_2$ (as defined in Definition 1.2) are given by
$$\text{Ran}(T_1)=E_1(A\ot A)
	\qquad\quad\text{and}\qquad\quad\
		\text{Ran}(T_2)=E_2(A\ot A)$$
where $E_1$ and $E_2$ are the idempotent maps from  $A\ot A$ to itself, given by left, respectively right multiplication by $E$. The element $E$ is uniquely determined by these conditions. We also assume that $E$ satisfies Assumption 1.10. The kernels of the canonical maps have the form
$$\text{Ker}(T_1)=(1 - G_1)(A\ot A)
	\qquad\quad\text{and}\qquad\quad\
		\text{Ker}(T_2)=(1-G_2)(A\ot A)$$
where $G_1$ and $G_2$ are the idempotent maps on $A\ot A$ determined by the equations
$$\align (G_1\ot \iota)(\Delta_{13}(a)(1\ot b\ot c))
	&= \Delta_{13}(a)(1\ot E)(1\ot b\ot c)\tag"(2.1)"\\
	(\iota\ot G_2)((a\ot b\ot 1)\Delta_{13}(c))
		&=(a\ot b\ot 1)(E\ot 1)\Delta_{13}(c)\tag"(2.2)"
\endalign$$
for all $a,b,c\in A$. See Proposition 1.11. Recall that we use $1$ here for the identity map.
\snl
Recall that in the regular case, it is expected that there is   
a right multiplier $F_1$ of $A\ot A^{\text{op}}$ and a left multiplier $F_2$ of $A^{\text{op}}\ot A$ so that the maps $G_1$ and $G_2$ are given as
$$G_1(a\ot b)=(a\ot 1)F_1(1\ot b)
\qquad\quad\text{and}\qquad\quad
G_2(a\ot b)=(a\ot 1)F_2(1\ot b)$$
for all $a,b\in A$. This will be investigated later in Section 4, see Proposition 4.5.
\nl
\it The generalized inverses $R_1$ and $R_2$ \rm
\nl
We will determine {\it generalized inverses} of the canonical maps, relative to these idempotent maps on the ranges and kernels. We will use the following standard result. Observe that we are using the symbols $r$, $t$ and $e$ for other things than usual in this paper. Also see e.g. [G] for more information about such generalized inverses.

\inspr{2.1} Lemma \rm 
Let $X$ be a vector space and $t$ a linear map from $X$ to itself. Assume that $e$ and $f$ are idempotent maps on $X$ so that $e$ projects on the range of $t$ and $1-f$ projects on the kernel of $t$. Then there is a unique linear map $r$ on from $X$ to itself so that
$$tr\xi=e\xi\qquad\quad\text{and} \qquad\quad rt\xi=f\xi$$
for all $\xi\in X$. Moreover we have  $trt=t$ and $rtr=r$.
\hfill$\square$\einspr

The proof is standard. The map $r$ is determined on the subspace $eX$ by the equation $rt\xi=f\xi$. And further it satisfies $r(1-e)=0$.
\snl
From the uniqueness, it follows that $r$ inherits properties from the maps $t$, $e$ and $f$ as we will see in the proof of the  application of this result to the canonical maps $T_1$ and $T_2$ with the given projection maps above (see Proposition 2.3 below). 
\nl
For the maps $T_1$ and $T_2$ we have the following properties. They behave with respect to multiplication as
$$T_1(a\ot bb')=(T_1(a\ot b))(1\ot b') 
	\qquad\quad\text{and} \qquad\quad
		T_2(a'a\ot b)=(a'\ot 1)T_2(a\ot b),$$
for $a,a',b,b'\in A$, and with respect to the coproduct as
$$(\Delta\ot\iota)T_1=(\iota\ot T_1)(\Delta\ot\iota)
	\qquad\quad\text{and} \qquad\quad
		(\iota\ot\Delta)T_2=(T_2\ot \iota)(\iota\ot\Delta).$$
In fact, the last two formulas are the same when they are written as 
$$(T_2\ot \iota)(\iota\ot T_1)=(\iota\ot T_1)(T_2\ot \iota)$$
and this is nothing else but coassociativity of the coproduct.
\snl
Now, we show that also the projection maps on the ranges and the kernels have the same behavior.

\inspr{2.2} Proposition \rm
For all $a,a',b,b'\in A$ we have
$$E_1(a\ot bb')=(E_1(a\ot b))(1\ot b') 
	\qquad\quad\text{and} \qquad\quad
		E_2(a'a\ot b)=(a'\ot 1)(E_2(a\ot b))$$
and the same for the maps $G_1$ and $G_2$. Furthermore
$$(\Delta\ot\iota)E_1=(\iota\ot E_1)(\Delta\ot\iota)
	\qquad\quad\text{and} \qquad\quad
		(\iota\ot\Delta)E_2=(E_2\ot \iota)(\iota\ot\Delta)$$
and the same for the maps $G_1$ and $G_2$.

\snl \bf Proof\rm:
Because $E_1$ is left multiplication with $E$ and $E_2$ is right multiplication with $E$, the first pair of formulas are obvious. It also follows immediately from the defining formulas (2.1) and (2.2) that $G_1$ and $G_2$ have the same behavior with respect to multiplication. See also a remark following the proof of Proposition 1.11.
\snl
Next, because $(\Delta\ot\iota)(E)=(1\ot E)(E\ot 1)$ we have for all $a,b\in A$ that
$$\align (\Delta\ot\iota)(E_1(a\ot b)) 
	&=(\Delta\ot\iota)(E(a\ot b))\\
	&=(\Delta\ot\iota)(E)(\Delta(a)\ot b)\\
	&=(1\ot E)(E\ot 1)(\Delta(a)\ot b)\\
	&=(1\ot E)(\Delta(a)\ot b)\\
	&=(\iota\ot E_1)((\Delta\ot\iota)(a\ot b)).
\endalign$$
This proves the first formula of the second pair. The second one uses also that \newline
$(\iota\ot\Delta)(E)=(1\ot E)(E\ot 1)$ and now $E_2(a\ot b)=(a\ot b)E$.
\snl
Finally, we show the equations with $G_1$ and $G_2$ (mentioned already in a remark following Proposition 1.11 in Section 1). We do it only for $G_1$ as the other case is completely similar. We have for all $a,b,c\in A$ that
$$\align (\Delta\ot\iota\ot\iota)(G_1\ot\iota)(\Delta_{13}(a)&(1\ot b\ot c))\\
	&=(\Delta\ot\iota\ot\iota)(\Delta_{13}(a)(1\ot E)(1\ot b\ot c))\\
	&=(\Delta\ot\iota\ot\iota)(\Delta_{13}(a))((1\ot 1\ot E)(1\ot 1\ot b\ot c))\\
	&=(\iota\ot\Delta_{13})(\Delta(a))((1\ot 1\ot E)(1\ot 1\ot b\ot c))\\
	&=(\iota\ot G_1\ot\iota)((\iota\ot\Delta_{13})(\Delta(a))(1\ot 1\ot b\ot c))\\
	&=(\iota\ot G_1\ot\iota)(\Delta\ot\iota\ot\iota)(\Delta_{13}(a)(1\ot b\ot c))
\endalign$$
and we see that $(\Delta\ot\iota)G_1=(\iota\ot G_1)(\Delta\ot\iota)$. This completes the proof.
\hfill$\square$\einspr

Now, we are ready to apply these results and get appropriate generalized inverses $R_1$ and $R_2$ of the canonical maps $T_1$ and $T_2$ respectively.  

\inspr{2.3} Proposition \rm
There are unique linear maps $R_1$ and $R_2$ so that
$$\align T_1R_1&=E_1 \qquad\qquad\text{and} \qquad\qquad R_1T_1=G_1 \\
T_2R_2&=E_2\qquad\qquad\text{and} \qquad\qquad R_2T_2=G_2.
\endalign$$
Moreover, also these (generalized) inverses satisfy, for all $a,a',b,b'\in A$,
$$\align R_1(a\ot bb')&=(R_1(a\ot b))(1\ot b') 
		\qquad\quad\text{and} \qquad\quad 
		(\Delta\ot\iota)R_1=(\iota\ot R_1)(\Delta\ot\iota)\\
	R_2(a'a\ot b)&=(a'\ot 1)(R_2(a\ot b))
		\qquad\quad\text{and} \qquad\quad
		(\iota\ot\Delta)R_2=(R_2\ot \iota)(\iota\ot\Delta).
\endalign$$

\snl\bf Proof\rm:
The existence of the generalized inverses $R_1$ and $R_2$ is a straightforward application of the general Lemma 2.1. 
\snl
To show that the inverses behave like the original maps with respect to multiplication and comultiplication we just apply the general principle formulated above. Consider e.g.\ the first formula. We have for $a,b,b'\in A$ that
$$\align R_1(a\ot bb')
	&=(R_1T_1R_1)(a\ot bb')\\
	&=R_1((T_1R_1(a\ot b))(1\ot b'))\\
	&=(R_1T_1)((R_1(a\ot b))(1\ot b'))\\
	&=(R_1T_1R_1(a\ot b))(1\ot b')\\
	&=(R_1(a\ot b))(1\ot b').
\endalign$$
We first have used the property for $T_1R_1$, then for $T_1$ and finally for $R_1T_1$.
\snl
Completely the same argument works for $R_2$. To prove the formulas with the coproduct, it is safer to look at the equivalent set of equations
$$\align (T_2\ot\iota)(\iota\ot R_1)&=(\iota\ot R_1)(T_2\ot\iota)\\
	(\iota\ot T_1)(R_2\ot\iota)&=(R_2\ot \iota)(\iota\ot T_1).
\endalign$$
Then the proof goes completely as before.
\hfill$\square$\einspr
\nl
\it The antipodes $S_1$ and $S_2$ \rm
\nl
We will now follow the treatment already given in [VD-W3, Proposition 2.4] to construct the antipodes $S_1$ and $S_2$ given the maps $R_1$ and $R_2$ as in the previous proposition. The correct interpretation of the formulas will be clear from the proof. We also give some more comments after the proof.

\inspr{2.4} Proposition \rm
There is a unique linear map $S_1$ from $A$ to $L(A)$, the left multipliers of $A$, and a unique linear map $S_2$ from $A$ to $R(A)$, the right multipliers of $A$, such that 
$$\align R_1(a\ot b)&=\sum_{(a)} a_{(1)} \ot S_1(a_{(2)})b \tag"(2.3)"\\
	R_2(a\ot b)&=\sum_{(b)}aS_2(b_{(1)})\ot b_{(2)} \tag"(2.4)"
\endalign$$
for all $a,b\in A$. 

\inspr{} Proof\rm:
We give the proof for $S_1$. The results for $S_2$ are proven in the same way.
\snl
Take $a,b\in A$ and define $S_1(a)b=(\varepsilon\ot\iota)(R_1(a\ot b))$. Because $$R_1(a\ot bb')=(R_1(a\ot b))(1\ot b')$$ for any other $b'\in A$ (cf.\ Proposition 2.3), we see that indeed $S(a)$ is a left multiplier of $A$ and this justifies the notation we have used.
\snl
Next, take $a,b,c\in A$. Then we have
$$\align \sum_{(a)} ca_{(1)}\ot S_1(a_{(2)})b 
	&=\sum_{(a)} (\iota\ot\varepsilon\ot\iota)(ca_{(1)}\ot R_1(a_{(2)}\ot b))\\
    &= (\iota\ot\varepsilon\ot\iota)(\iota\ot R_1)(T_2\ot\iota)(c\ot a\ot b)).
\endalign$$
Now we use that $(\iota\ot R_1)(T_2\ot\iota)=(T_2\ot\iota)(\iota\ot R_1)$ (cf.\ Proposition 2.3) and we obtain
$$\sum_{(a)} ca_{(1)}\ot S_1(a_{(2)})b 
=(\iota\ot\varepsilon\ot\iota)(T_2\ot\iota)(\iota\ot R_1)(c\ot a\ot b)).$$
Because $(\iota\ot\varepsilon)(T_2(p\ot q))=pq$ for all $p,q$, this last expression is equal to \newline
$(c\ot 1)R_1(a\ot b)$. This proves formula (2.3) of the proposition.
\snl
Moreover, we see from the proof that $S_1$ is completely determined by the given property (and so it is unique).
\hfill$\square$
\einspr

The following are important remarks.

\inspr{2.5} Remark \rm
i) Whereas we have $R_1(a\ot b)$ belongs to $A\ot A$, this is not obvious for the right hand side $\sum_{(a)} a_{(1)} \ot S_1(a_{(2)})b$ in the equation (2.3). However, it follows from the argument in the proof that this is the case. 
\snl
ii) We will need this for the interpretation of the formulas in the next proposition. Consider the first formula in the next proposition with $S=S_1$. Take $b\in A$. We have that $\Delta(a)(1\ot b)$ is in $A\ot A$. Write this element as $\sum_i p_i \ot q_i$ with $p_i,q_i\in A$. Then we have
$$\sum_{(a)} a_{(1)}S_1(a_{(2)})a_{(3)}b=\sum_{i,(p_i)}p_{i(1)}S_1(p_{i(2)})q_i$$
and this is well-defined in $A$ as $\sum_{(p)}p_{(1)}\ot S_1(p_{(2)})q$ is in $A\ot A$ for all $p,q\in A$.
\snl
iii) For the second formula below (in Proposition 2.6), again with $S=S_1$, take once more $b\in A$ and use first that $\sum_{(a)}a_{(1)}\ot S_1(a_{(2)})b$ is in $A\ot A$. If we now write this element as $\sum_i p_i \ot q_i$ with $p_i,q_i\in A$, we find
$$\sum_{(a)} S_1(a_{(1)})a_{(2)}S_1(a_{(3)})b=\sum_{i,(p_i)}S_1(p_{i(1)})p_{i(2)}q_i$$
and again this is well-defined in $A$.
\hfill$\square$\einspr

The same phenomenon occurs in the theory of multiplier Hopf algebras. However, as we are usually working with regular multiplier Hopf algebras, so that the antipode $S$ is bijective, we then write e.g.
$$ \sum_{(a)} a_{(1)} \ot S(a_{(2)})b=(\iota \ot S)((1\ot c)\Delta(a))$$ where $c=S^{-1}(b)$ and then it is clear that this element is in $A\ot A$. This is one of the peculiarities that can happen in the non-regular case.
\nl
Now, the following two well-known formulas in the theory of Hopf algebras, are also true here, but they need the correct interpretation (as explained in the remark above). 

\inspr{2.6} Proposition \rm
The following two equations
$$\sum_{(a)} a_{(1)}S(a_{(2)})a_{(3)}=a
\qquad\quad\text{and}\quad\qquad
         \sum_{(a)} S(a_{(1)})a_{(2)}S(a_{(3)})=S(a)
$$
for all $a\in A$ are satisfied for $S$ equal to $S_1$ and $S_2$. In the case of $S_1$, we have two equalities as left multipliers (see the remark above) whereas in the case of $S_2$ we have equalities as right multipliers.

\snl\bf Proof\rm:
The two equations for $S_1$ follow if we insert the formula for $R_1$ in terms of $S_1$, given in the previous proposition, in the formulas $T_1R_1T_1=T_1$ and $R_1T_1R_1=R_1$, and if we then apply $\varepsilon$ (or the fullness of $\Delta$). Similarly for $S_2$.
\hfill $\square$
\einspr

The formulas can also be interpreted as $\iota * S * \iota= \iota$ and $S * \iota * S= S$ in the convolution algebra of linear maps from $A$ to itself, but of course, it would be more difficult to make these formulas meaningful. Nevertheless, it is good to keep this interpretation in mind.
\snl
Needless to say that  maps $S_1$ and $S_2$, satisfying the formulas of the previous proposition, yield  generalized inverses $R_1$ and $R_2$ of
$T_1$ and $T_2$ satisfying the conditions of Proposition 2.3.
\nl
\it The antipode $S:A\to M(A)$\rm
\nl
We will now show that the two antipodes $S_1$ and $S_2$ we obtained in the previous item {\it coincide}. It follows from the {\it specific choices} of the maps $G_1$ and $G_2$ as we will see in the proof. As the following result is very crucial, we will give more remarks after we have proven it.

\inspr{2.7} Proposition \rm
The antipodes $S_1$ and $S_2$ satisfy $b(S_1(a)c)=(bS_2(a))c$ for all $a,b,c\in A$ and hence define a map $S:A\to M(A)$,  called {\it the antipode} of the weak multiplier Hopf algebra $(A,\Delta)$.

\snl\bf Proof\rm:
Take $a\in A$. For all $p,q\in A$ we have
$$\align
	 \sum_{(a)} a_{(1)}\ot a_{(2)} S_1(a_{(3)})p \ot  a_{(4)}q 
	&= (T_1R_1\ot \iota) (\sum_{(a)} a_{(1)}\ot p \ot  a_{(2)}q). \\
	&= (E\ot 1) \Delta_{13}(a)(1\ot p\ot q). 
\endalign$$
On the other hand, we have for all $r,s\in A$ that
$$\align
	 \sum_{(a)} ra_{(1)} \ot sa_{(2)} S_2(a_{(3)})\ot  a_{(4)} 
	&= (\iota\ot R_2T_2)(\sum_{(a)} ra_{(1)} \ot s \ot  a_{(2)}) \\
	&= (\iota\ot G_2) ((r\ot s\ot 1)\Delta_{13}(a)). 
\endalign$$
By the definition of $G_2$, this last expression is
$$(r\ot s\ot 1)(E\ot 1)\Delta_{13}(a)$$
and if we combine the two results we get
$$ \sum_{(a)} ra_{(1)}\ot (sa_{(2)} S_1(a_{(3)}))p \ot  a_{(4)}q=
\sum_{(a)} ra_{(1)}\ot s(a_{(2)} S_2(a_{(3)})p) \ot  a_{(4)}q
$$
for all $p,q,r,s\in A$. It follows that 
$$\sum_{(a)} (sa_{(1)} S_2(a_{(2)}))p =
 \sum_{(a)} s(a_{(1)} S_1(a_{(2)})p)$$
for all $p,s\in A$. 
\snl
In a completely similar way, using the definition of the map $G_1$ we find
$$ \sum_{(a)} s(S_1(a_{(1)}) a_{(2)}p) =
\sum_{(a)} (sS_2(a_{(1)})a_{(2)})p 
$$
for all $p,s\in A$. 
\snl
Then, we get for all $b,c\in A$ that
$$\align 
	b(S_1(a)c) &= \sum_{(a)} b(S_1(a_{(1)}) a_{(2)} S_1(a_{(3)})c) \\ 
		&= \sum_{(a)} (bS_2(a_{(1)}) a_{(2)}) S_1(a_{(3)})c \\
		&= \sum_{(a)} (bS_2(a_{(1)})) (a_{(2)} S_1(a_{(3)})c) \\
		&= \sum_{(a)} (bS_2(a_{(1)}) a_{(2)} S_2(a_{(3)}))c \\
		&= (bS_2(a))c.
\endalign$$
This proves the result.
\hfill$\square$\einspr

Because this result is very important for our treatment, and in some sense also quite subtle, we add the following remarks.

\inspr{2.8} Remark \rm
i) First of all, we have used the Sweedler notation in the proof, but we have made sure that everything is well-covered. In fact, it would be possible (and straightforward) to rewrite the proof in terms of the maps $T_1, T_2, R_1$ and $R_2$, without the Sweedler notation. Of course this would be much less transparent. We refer also to the remarks made in 2.5.
\snl
ii) As we see from the last argument in the proof, the equality $S_1=S_2$ follows from the equations (and hence is equivalent with)
$$\sum_{(a)} a_{(1)} S_1(a_{(2)})= \sum_{(a)} a_{(1)} S_2(a_{(2)})
\quad\quad\text{and}\quad\quad 
\sum_{(a)} S_1(a_{(1)}) a_{(2)}= \sum_{(a)} S_2(a_{(1)}) a_{(2)}$$
for all $a$. 
\snl
iii) These equalities in turn are equivalent with the defining formulas (2.1) and (2.2) of the maps $G_1$ and $G_2$. We refer also to the discussion about this topic in [VD-W3] (see Proposition 3.17 in that paper and the Remark 3.18 following it).
\snl
iv) In other words, the equality $S_1=S_2$ determines the kernels of the canonical maps $T_1$ and $T_2$, given their ranges. We feel that this is {\it quite remarkable} (and crucial for the theory).
\hfill$\square$\einspr

We have now constructed the antipode and proven its first  properties. More (and important) results about the antipode and some consequences of this will be postponed till the next section. Now, we see how a weak multiplier Hopf algebra can be characterized  when an antipode is given. This will be used later in this section to show that any weak Hopf algebra (as introduced and studied in [B-N-S], see also [N-V1]) is a weak multiplier Hopf algebra (as defined in Section 1).
\nl
\it A characterization in terms of the antipode \rm
\nl
Before we look at our examples, we first obtain a characterization of weak multiplier Hopf algebras in terms of an existing antipode. It looks more like the definition of a weak Hopf algebra. Moreover, it is often more convenient to show in a concrete case that we have a weak multiplier Hopf algebra as most of the time, the antipode is already available.
\snl
We first state the result and prove it. Then we will give some more comments.

\inspr{2.9} Theorem \rm
Let $(A,\Delta)$ be a pair of a non-degenerate idempotent algebra $A$ with a full coproduct $\Delta$ and such that there is a counit. 
\snl
i) Assume that there is a linear map $S:A\to M(A)$ such that the maps $R_1$ and $R_2$, defined by the formulas (2.3) and (2.4) in Proposition 2.4, have range in $A\ot A$ and so that, for all $a\in A$, 
$$\sum_{(a)} a_{(1)}S(a_{(2)})a_{(3)}=a
\qquad\quad\text{and}\quad\qquad
         \sum_{(a)} S(a_{(1)})a_{(2)}S(a_{(3)})=S(a)\tag"(2.5)"
$$
hold in $M(A)$ (as in Proposition 2.6).
\snl
ii) Assume that there is an element $E$ in $M(A\ot A)$ so that
$$(T_1R_1)(a\ot b)=E(a\ot b) 
\qquad\quad\text{and}\quad\qquad
(T_2R_2)(a\ot b)=(a\ot b)E\tag"(2.6)"$$
for all $a,b\in A$ (as in Proposition 2.3) and that 
$$(\iota\ot \Delta)(E)=(\Delta\ot\iota)(E)=(E\ot 1)(1\ot E)=(1\ot E)(E\ot 1)\tag"(2.7)"$$
(as in condition ii) of Definition 1.14).
\snl
Then $(A,\Delta)$ is a weak multiplier Hopf algebra and $S$ is its antipode.

\snl\bf Proof\rm:
First remark that the formulas (2.5) make sense by the requirement that the maps $R_1$ and $R_2$ have range in $A\ot A$ (cf.\ Remark 2.5). Also remember that these formulas essentially are equivalent with the fact that $R_1$ and $R_2$ are generalized inverses of $T_1$ and $T_2$ respectively (cf.\ a remark following Proposition 2.6).
\snl
Then it follows from the formulas (2.6) in condition ii) that $E$ is an idempotent in $M(A\ot A)$ and that 
$$E(A\ot A)=T_1(A\ot A)
	\qquad\quad\text{and}\qquad\quad
		(A\ot A)E=T_2(A\ot A).$$
We see that already the conditions i) and ii) of Definition 1.14 are fulfilled.
\snl
It remains to show that finally, also condition iii) of Definition 1.14 is true.
\snl
To prove this, first remark that the kernels of the maps $T_1$ and $T_2$ are given by the ranges of the idempotents $1-R_1T_1$ and $1-R_2T_2$ because $R_1$ and $R_2$ are generalized inverses of $T_1$ and $T_2$ respectively. Therefore, we only need to show that the idempotent maps $R_1T_1$ and $R_2T_2$ satisfy the formulas for $G_1$ and $G_2$ as given in Proposition 1.11, see also (2.1) and (2.2).
\snl
We will prove this for the first map. The proof for the second one is completely similar. Also compare with the argument given in the proof of Proposition 2.7.
\snl
Take $a,b,c\in A$. Then we have
$$\align
	(R_1T_1\ot \iota)(\Delta_{13}(a)(1\ot b\ot c))
	&=\sum_{(a)}R_1T_1 (a_{(1)}\ot b) \ot a_{(2)}c \\
	&=\sum_{(a)} a_{(1)}\ot S(a_{(2)}) a_{(3)}b \ot a_{(4)}c. 
\endalign$$
If we multiply from the left with $p\ot q\ot 1$, where $p,q\in A$, we can write the right hand side as 
$$\align 
	\sum_{(a)} pa_{(1)}\ot (T_2R_2(q\ot a_{(2)}))(b\ot c)
 	&=\sum_{(a)} pa_{(1)}\ot (q\ot a_{(2)})E(b\ot c)\\
	&=(p\ot q\ot 1)\Delta_{13}(a)(1\ot E)(1\ot b\ot c).
\endalign$$
If we insert this expression again in the first series of equations, and cancel $p,q$, we precisely obtain the required formula. 
\snl
Finally, it is clear that $S$ is the antipode of the weak multiplier Hopf algebra $(A,\Delta)$.
\hfill$\square$\einspr

We will use this result to show that any weak Hopf algebra is a weak multiplier Hopf algebra. Indeed, the above characterization is already very close to the defining properties for a weak Hopf algebra.
\snl
If we compare the conditions in the proposition with those in the original definition of a weak multiplier Hopf algebra, we see that the existence of the antipode is used to prove the assumption iii) in Definition 1.14 about the kernels of the canonical maps. The other conditions i) and ii) of the original definition appear essentially also in the new characterization. It means that the original axioms are, in a sense, weaker than the ones in the characterization of the above theorem.
\snl
On the other hand, it is precisely condition iii) of the original definition that is not so obvious to show in concrete examples and in fact, the existing antipode helps to do this job. This means that the original definition is nice from a theoretical point of view whereas the characterization of the previous theorem is nicer from a practical point of view. A similar statement is true already for multiplier Hopf algebras.
\snl
Finally remark once more how the conditions on the kernels of the canonical maps are intimately related with the existence of the antipode as we see this (again) at the end of the proof above.
\snl
This seems to be the right place to refer to a first, unpublished and incomplete version of an attempt to define and develop this theory of weak multiplier Hopf algebras. Indeed, in [VD-W2], a notion of {\it Multiplier Unifying Hopf algebras} has been considered, with axioms like (2.5) in Theorem 2.9.
\nl
\it Special cases and examples \rm
\nl
We begin with the weak Hopf algebras. 
We formulate the result as a proposition.

\iinspr{2.10} Proposition \rm
Let $(A,\Delta)$ be a weak Hopf algebra. Then it is a weak multiplier Hopf algebra. 

\bf\snl Proof\rm:
Assume that $(A,\Delta)$ is a weak Hopf algebra.
\snl
Because $A$ is an algebra with identity by assumption, it is non-degenerate and idempotent. The coproduct $\Delta$ satisfies our requirements (as in Definition 1.1). Furthermore, by assumption there is a counit and in this case, it follows that the coproduct is also full. The coproduct is also automatically regular in this case.
\snl
By definition, there exists an antipode $S$ and it satisfies the conditions i) of Theorem 2.9. Indeed, $S$ is a linear map from $A$ to $A$ and of course in this case, the associated maps $R_1$ and $R_2$ have range in $A\ot A$. The second formula in (2.5) is  part of the axioms of a weak Hopf algebra. The first one is not, but easy to obtain from the axioms. Indeed, we have for all $a\in A$ that
$$\align
\sum_{(a)}a_{(1)}S(a_{(2)})a_{3)}
&=\sum_{(a)}(\varepsilon\ot\iota)(\Delta(1)(a_{(1)}\ot 1))a_{(2)} \\
&=\sum_{(a)}(\varepsilon\ot\iota)(\Delta(1)(a_{(1)}\ot a_{(2)}) \\
&=(\varepsilon\ot\iota)\Delta(a)=a.
\endalign$$
This argument can be found e.g.\ already in [B-N-S]. 
\snl
For $E$ we take of course $\Delta(1)$. Again by definition, the formula (2.7) is fulfilled. So, we only have to verify the formulas in (2.6). Now we have for $a\in A$ that
$$\align T_1R_1(a\ot 1)
&=\sum_{(a)} a_{(1)}\ot  a_{(2)} S(a_{(3)})\\
&=\sum_{(a)} a_{(1)}\ot  (\varepsilon\ot\iota)(\Delta(1)(a_{(2)}\ot 1))\\
&=(\iota\ot\varepsilon\ot\iota)((1\ot\Delta(1))(\Delta(a)\ot 1))\\
&=(\iota\ot\varepsilon\ot\iota)((1\ot\Delta(1))(\Delta(1)\ot 1)(\Delta(a)\ot 1))\\
&=(\iota\ot\varepsilon\ot\iota)((\Delta\ot\iota)(\Delta(1)))(\Delta(a)\ot 1))\\
&=(\iota\ot\varepsilon\ot\iota)((\Delta\ot\iota)(\Delta(1)(a\ot 1)))\\
&=\Delta(1)(a\ot 1).
\endalign$$
This shows that $T_1R_1(a\ot b)=E(a\ot b)$. This argument is also found in the original paper [B-N-S].
\snl
Similarly we get $T_2R_2(a\ot b)=(a\ot b)E$. 
\hfill$\square$\einspr

Remark that originally, weak Hopf algebras were only considered in the finite-dimensional case (see e.g.\ [B-N-S] and also [N-V2]).  Later however, also infinite-dimensional weak Hopf algebras have been studied (see e.g.\ [N]). Indeed, it turns out that the axioms, as first given in [B-N-S] and many of the results do not depend on the finite-dimensionality of the underlying algebra (as we see above). 
\snl
It is known that the antipode of a finite-dimensional weak Hopf algebra is bijective (see e.g.\ Proposition 2.10 in [B-N-S]). As we will show in the Section 4, this means that we get a {\it regular} weak multiplier Hopf algebra. Conversely, we will prove that any regular weak multiplier Hopf algebra with a finite-dimensional underlying algebra must be a weak Hopf algebra (see Proposition 4.12  in Section 4). In fact, we will show that any regular weak multiplier Hopf algebra with a unital underlying algebra is a weak Hopf algebra.
\snl
However, we were not able to show that any weak multiplier Hopf algebra with a finite-dimensional underlying algebra is automatically regular. More generally, we do not know if any weak multiplier Hopf algebra with an underlying untial algebra is a weak Hopf algebra. This raises the question whether or not, our axioms are weaker. 
\snl
We will discuss this further at the end of the Section 4. 
\nl
Now we have a brief look at the examples coming from a groupoid. We have mentioned in Section 1 already that they are examples of weak multiplier Hopf algebras, but we did not give a proof in detail. What we do here is just illustrating how the characterization in Theorem 2.9 can be used to show this.

\iinspr{2.11} Example \rm
Take a groupoid $G$ and the algebra $K(G)$ of complex functions on $G$ with finite support and pointwise product as in Example 1.15. The antipode is given by $(S(f))(p)=f(p^{-1})$ for all $f\in K(G)$ and $p\in G$. We will verify the required conditions in Theorem 2.9.
\snl
Because $p=pp^{-1}p$ and of course also $p^{-1}=p^{-1}pp^{-1}$ for all elements $p$ in $G$, the antipode will satisfy the formulas (2.5) in Theorem 2.9. 
\snl
Moreover, we find that, for all $f,g\in K(G)$ and all $p,q\in G$,
$$(\sum_{(f)} f_{(1)} \ot f_{(2)} S(f_{(3)})g)(p,q)
=\sum_{(f)} f_{(1)}(p) (f_{(2)} S(f_{(3)})g)(q)
=f(pqq^{-1})g(q).$$
We see that this is $f(p)g(q)$ if $pq$ is defined and $0$ otherwise. Recall that in Example 1.15, we defined $E$ as multiplication with the function on pairs $(p,q)$ having the value $1$ if $pq$ is defined and $0$ otherwise. So we find
$$\sum_{(f)} f_{(1)} \ot f_{(2)} S(f_{(3)})g=E(f\ot g)$$
for all $f,g\in K(G)$. 
Similarly, we find
$$(\sum_{(f)} f_{(1)} \ot S(f_{(2)}) f_{(3)}g)(p,q)
=f(pq^{-1}q)g(q)$$
for all $f,g\in K(G)$ and $p,q\in G$. Then, for all $f,g$,
$$\sum_{(f)} f_{(1)} \ot S(f_{(2)}) f_{(3)}g=F_1(f\ot g)$$
where $F_1$ is multiplying with the function on pairs $(p,q)$ that takes the value $1$ if $s(p)=s(q)$ and $0$ on other pairs. Also this is in agreement with the results in Example 1.15.
Similar calculations give, for all $f,g\in K(G)$,
$$\align 
\sum_{(g)} fS(g_{(1)})g_{(2)}\ot g_{(3)}&=(f\ot g)E\\
\sum_{(g)} fg_{(1)}S(g_{(2)})\ot g_{(3)}&=F_2(f\ot g)
\endalign$$
where now $F_2$ is the function on pairs $(p,q)$ that is $1$ if $t(p)=t(q)$ and $0$ otherwise. Again, these two formulas are in agreement with the definitions and results in Example 1.15.
\snl
Certainly, we have shown that all the conditions of Theorem 2.9 are fulfilled (remark that the algebra is abelian here) and therefore that we get a weak multiplier Hopf algebra.
\hfill$\square$\einspr

The dual case, as constructed in Example 1.16, is treated in a completely similar way. The antipode on the convolution algebra $\Bbb C G$ is now given by $S(\lambda_p)=\lambda_{p^{-1}}$ for all $p\in G$ and where $p\mapsto \lambda_p$ is the canonical imbedding of $G$ in $\Bbb C G$. There is a slight complication as the algebra $\Bbb C G$ might be non-abelian, contrary to the algebra $K(G)$ above. On the other hand, the legs of $E$ in this case are an abelian subalgebra of the multiplier algebra and this again simplifies the arguments.

\nl\nl



\bf 3. The antipode: main results\rm
\nl
We start again with a weak multiplier Hopf algebra as defined in Definition 1.14. We consider the antipode $S: A\to M(A)$ as obtained in the previous section. We have already shown some elementary properties. Now, in this section, we will prove the main results about this antipode. In particular, we show that it is a non-degenerate anti-homomorphism. Then it has a unique extension to a unital anti-homomorphism from $M(A)$ to itself. We also show that it is an anti-coalgebra map (in a sense to be made precise). To prove these results, we need to introduce the source and target maps $\varepsilon_s$ and $\varepsilon_t$ and some properties of their images in $M(A)$. More results about these objects are obtained in a separate paper (see [VD-W4]). 
\nl
We know already from previous remarks that 
$$\sum_{(a)}S(a_{(1)})a_{(2)}
\qquad\quad\text{and}\qquad\quad
\sum_{(a)}a_{(1)}S(a_{(2)})
$$
are well-defined multipliers in $M(A)$ for any $a\in A$ (also see Remark 2.5 in the previous section). Because we will use these expressions at various places in the forthcoming proofs, let us introduce the following notations. 

\inspr{3.1} Definition \rm 
For $a$ in $A$ we set
$$\varepsilon_s(a)=\sum_{(a)}S(a_{(1)})a_{(2)}
\qquad\quad\text{and}\quad\qquad
\varepsilon_t(a)=\sum_{(a)}a_{(1)}S(a_{(2)}).$$
The map $\varepsilon_s$ is called the {\it source map} while $\varepsilon_t$ is the {\it target map}.
\hfill$\square$\einspr 

We will use these maps further for proving the main results about the antipode, but we will not need all there nice properties in this paper. These will be considered in [VD-W4]. We just collect here what we need to obtain the results about the antipode.
\snl

First we observe in the following lemma that the images $\varepsilon_s(A)$ and $\varepsilon_t(A)$ are, in a way, the {\it left} and the {\it right} leg of $E$.

\inspr{3.2} Lemma \rm
We have that 
$$
(\omega(c\,\cdot\,a)\ot\iota)(E)\in \varepsilon_t(A)
\qquad\quad\text{and}\quad\qquad
(\iota\ot\omega(a\,\cdot\,c))(E)\in \varepsilon_s(A)
$$
for all linear functionals $\omega$ on $A$ and all $a,c\in A$. Moreover, any element in $\varepsilon_s(A)$ and $\varepsilon_t(A)$ respectively is obtained as a linear combination of such elements.

\bf\snl Proof\rm: We have 
$$E(a\ot b)=T_1R_1(a\ot b)=\sum_{(a)}a_{(1)}\ot a_{(2)}S(a_{(3)})b$$
for all $a,b\in A$ and so we can write
$$(c\ot 1)E(a\ot 1)=\sum_{(a)}ca_{(1)}\ot a_{(2)}S(a_{(3)})$$
for all $a,c\in A$. This element belongs to $A\ot \varepsilon_t(A)$. If we apply a linear functional $\omega$ on the first leg, we find the first property in the lemma. Similarly, the second one can be proven.
\snl
The second statement in the lemma follows from the fullness of the coproduct. 
\hfill$\square$\einspr


Because we have $E(a\ot 1)=\sum_{(a)}a_{(1)}\ot a_{(2)}S(a_{(3)})$ for all $A$, we also have 
$$E(aa'\ot 1)=\sum_{(a)}a_{(1)}a'\ot a_{(2)}S(a_{(3)})$$
for all $a,a'\in A$. Therefore, if the coproduct is {\it regular}, we see that already \newline 
$E(A\ot 1)\subseteq A\ot \varepsilon_t(A)$. We use that $A$ is idempotent. Similarly, in that case, also $(1\ot A)E\subseteq \varepsilon_s(A)\ot A$. This observation is important for the further study of the images of the source and target maps in [VD-W4.]

\inspr{3.3} Lemma \rm
For all $a\in A$, we have
$$\align \Delta(\varepsilon_t(a))&=E(\varepsilon_t(a)\ot 1)=(\varepsilon_t(a)\ot 1)E \tag"(3.1)"\\
\Delta(\varepsilon_s(a))&=E(1\ot \varepsilon_s(a))=(1\ot \varepsilon_s(a))E \tag"(3.2)". 
\endalign$$

\snl\bf Proof\rm: 
These formulas will follow from the equations
$$\align
(\iota\ot \Delta)E&=(1\ot E)(E\ot 1)=(E\ot 1)(1\ot E) \\
(\Delta\ot \iota)E&=(1\ot E)(E\ot 1)=(E\ot 1)(1\ot E),
\endalign$$
together with the fact that $\varepsilon_s(A)$ and $\varepsilon_t(A)$ are respectively the left and the right leg of $E$ in the sense explained in the proof of Lemma 3.2.
\hfill$\square$\einspr

The following is an easy consequence.

\inspr{3.4} Lemma \rm
The images $\varepsilon_s(A)$ and $\varepsilon_t(A)$ are commuting subalgebras of $M(A)$.

\snl\bf Proof\rm:
Take any element $a\in A$ and $y\in \varepsilon_s(A)$. Then we have $\Delta(ay)=\Delta(a)(1\ot y)$ because of the previous lemma. If we apply $S\ot \iota$ and multiply, we find that $\varepsilon_s(ay)=\varepsilon_s(a)y$. This shows that $\varepsilon_s(A)$ is a subalgebra of $M(A)$. Similarly for $\varepsilon_t(A)$. As the legs of $E$ commute, it follows that these algebras commute with each other.
\hfill$\square$\einspr

Now we are ready to prove the main properties of the antipode. First we have the following.

\inspr{3.5} Proposition \rm
We have $S(ab)=S(b)S(a)$ for all $a,b\in A$.

\snl\bf Proof\rm:
For any $a,b\in A$ we have
$$E(a\ot b)=\sum_{(a)}\Delta(a_{(1)})(1\ot S(a_{(2)})b).$$
Then we have for all $a,a',b\in A$ that on the one hand
$$E(aa'\ot b)= \sum_{(a)(a')}\Delta(a_{(1)}a'_{(1)})(1\ot S(a_{(2)}a'_{(2)})b)
$$
while on the other hand
$$\align E(aa'\ot b)&=E(a\ot b)(a'\ot 1)\\
&=\sum_{(a)}\Delta(a_{(1)})(a'\ot S(a_{(2)})b)\\
&=\sum_{(a)}\Delta(a_{(1)})E(a'\ot S(a_{(2)})b)\\
&=\sum_{(a)(a')}\Delta(a_{(1)})\Delta(a'_{(1)})(1\ot S(a'_{(2)})S(a_{(2)})b).
\endalign$$
If we apply $\varepsilon\ot\iota$ on both expressions and then cancel $b$, we get
$$\sum_{(a)(a')}a_{(1)}a'_{(1)}S(a_{(2)}a'_{(2)})=
\sum_{(a)(a')}a_{(1)}a'_{(1)}S(a'_{(2)})S(a_{(2)})$$
for all $a,a'\in A$. Similarly, we can obtain 
$$\sum_{(a)(a')}S(a_{(1)}a'_{(1)})a_{(2)}a'_{(2)}=
\sum_{(a)(a')}S(a'_{(1)})S(a_{(1)})a_{(2)}a'_{(2)}$$
for all $a,a'\in A$. If we now combine these results with the formulas in Proposition 2.6 in the previous section and with the result of the previous lemma, we get
$$\align S(aa')
&=\sum_{(a)(a')}S(a_{(1)}a'_{(1)})a_{(2)}a'_{(2)}S(a_{(3)}a'_{(3)})\\
&=\sum_{(a)(a')}S(a'_{(1)})S(a_{(1)})a_{(2)}a'_{(2)}S(a'_{(3)})S(a_{(3)})\\
&=\sum_{(a)(a')}S(a'_{(1)})a'_{(2)}S(a'_{(3)})S(a_{(1)})a_{(2)}S(a_{(3)})\\
&=S(a')S(a)
\endalign$$
for all $a,a'\in A$. This proves the result.
\hfill$\square$\einspr

We did not bother too much about coverings in the above proof. But in view of Remark 2.5.i, things will be all right if we multiply e.g.\ the last set of equalities with an element $b$ from the right and an element $c$ from the left.
\snl
So we see that, as expected, the antipode is an {\it anti-algebra map}. We will show that (in some sense) it is also an {\it anti-coalgebra map}. Before we do this however, we consider some more properties of the antipode as an anti-homomorphism from $A$ to $M(A)$.

\inspr{3.6} Proposition \rm
In any weak multiplier Hopf algebra, we have $AS(A)=A$ as well as $S(A)A=A$.

\snl \bf Proof\rm:
We know that 
$$\Delta(x)(a\ot b)=\sum_{(a)} \Delta(xa_{(1)})(1\ot  S(a_{(2)})b)$$
for all $a,b,x\in A$. We multiply with an element $c$ from the left in the first factor and cancel $b$. Then we find
$$(c\ot 1)\Delta(x)(a\ot 1)=\sum_{(a)} (c\ot 1)\Delta(xa_{(1)})(1\ot  S(a_{(2)})) \tag"(3.3)"$$
and this equation is true in $A\ot A$ for all $a,c,x\in A$. The expression on the right belongs to $A\ot AS(A)$. Now suppose that $\omega$ is a linear functional on $A$ that is $0$ on $AS(A)$. We will show that it has to be zero on all of $A$ and this will prove that $AS(A)=A$.
\snl
Indeed, apply $\iota\ot \omega$ to the equation (3.3). Then we get 
$$(\iota\ot \omega)((c\ot 1)\Delta(x)(a\ot 1))=0$$
for all $a,c,x\in A$. By the the non-degeneracy of the product in $A$, we can cancel $a$ and we find that also $(\iota\ot \omega)((c\ot 1)\Delta(x))=0$ for all $c,x\in A$. Then it follows from the fullness of $\Delta$ that $\omega=0$. This completes the proof of the statement $AS(A)=A$.
\snl
The other statement $S(A)A=A$ is proven in a similar way, now starting with the formula
$$(b\ot a)\Delta(x)=\sum_{(a)} (bS(a_{(1)})\ot 1) \Delta(a_{(2)}x)$$
for all $a,b,x\in A$.
\hfill$\square$\einspr

This result means that $S:A\to M(A)$ is a {\it non-degenerate} anti-homomorphism and by the general theory, it has a unique extension, still denoted by $S$ as usual, to a unital anti-homomorphism from $M(A)$ to itself. 
\snl
It is also a consequence of this result that the tensor product map $S\ot S:A\ot A\to M(A\ot A)$ is non-degenerate and so it can also be extended to a unital anti-homomorphism from $M(A\ot A)$ to itself. This result will be used below when we show that $S$ is an anti-coalgebra map as well.
\snl
Before we do this, let us make one more related remark. From the result above, it follows that the product in $S(A)$ is non-degenerate. Then we can consider the multiplier algebra $M(S(A))$. Again using this result, we see that the imbedding of $S(A)$ in $M(A)$ is non-degenerate and so extends to $M(S(A))$. It turns out that this is still an imbedding. This means that we can consider elements of $M(S(A))$ as sitting in $M(A)$. We will reconsider this in Section 2 of [VD-W4].
\nl
Next we will show that the antipode is an anti-coalgebra map.  Unfortunately, the result is not as nice as it is hoped for. We discuss the problem later (cf.\ Remark 3.8.iii further), as well as in the next section (see a remark follwong Proposition 4.6). 

\inspr{3.7} Proposition \rm For all $a\in A$, we have 
$$\Delta(S(a))=E(\sigma(S\ot S)\Delta(a))=(\sigma(S\ot S)\Delta(a))E.$$

\snl\bf Proof\rm:
Take $a\in A$. Then we have
$$\align 
	\Delta(S(a))
	&=\sum_{(a)}\Delta(S(a_{(1)})a_{(2)}S(a_{(3)}))\\
    &=\sum_{(a)}\Delta(S(a_{(1)}))\Delta(a_{(2)}S(a_{(3)}))\\
    &=\sum_{(a)}\Delta(S(a_{(1)}))E(a_{(2)}S(a_{(3)})\ot 1)\\
    &=\sum_{(a)}\Delta(S(a_{(1)}))\Delta(a_{(2)})(S(a_{(4)})\ot S(a_{(3)}))
\endalign$$
where we have used formula (3.1) from Lemma 3.3 and in the last equality that \newline
$E(p\ot 1)=\sum_{(p)}\Delta(p_{(1)})(1\ot S(p_{(2)}))$ with $a_{(2)}$ in the place of $p$. Then
$$\align 
	\Delta(S(a))
    &=\sum_{(a)}\Delta(S(a_{(1)})a_{(2)})(S(a_{(4)})\ot S(a_{(3)}))\\
    &=\sum_{(a)}E(1\ot(S(a_{(1)})a_{(2)}))(S(a_{(4)})\ot S(a_{(3)}))
\endalign$$
where now we have used formula (3.2) from Lemma 3.3. Finally, we get
$$\align 
	\Delta(S(a))
    &=\sum_{(a)}E(S(a_{(4)})\ot S(a_{(1)})a_{(2)}S(a_{(3)}))\\
    &=\sum_{(a)}E(S(a_{(2)})\ot S(a_{(1)})).
\endalign$$
In a similar way, we find such a formula with $E$ on the other side.
\hfill$\square$\einspr

Of course, to give a meaning to all these equations, we need to multiply at the right places with the right elements in order to get a good covering of the expressions. However, it is not hard to do this, taking into account all earlier considerations about this problem.
\snl
Before we continue, we want to make a few more remarks on this result. 

\inspr{3.8} Remark \rm
i) Obviously, the left hand side of the equation in Proposition 3.7 presents no problem because $S(a)\in M(A)$ and we have extended the coproduct to the multiplier algebra. Similarly, the right hand side is fine because we also have extended $S\ot S$ to $M(A\ot A)$. 
\snl
ii) There is another way to give a meaning to the right hand side in  $M(A\ot A)$. To do this, remark that the two expressions
$$\sum_{(a)} a_{(1)} \ot S(a_{(2)})b
  \qquad\quad\text{and}\qquad\quad
    \sum_{(a)} cS(a_{(1)})\ot a_{(2)}$$
are well-defined in $A\ot A$ for all $a,b,c\in A$ (see Remark 2.5 in the previous section). Then, because the antipode $S$ maps into the multiplier algebra, we also get that the two expressions
$$ \sum_{(a)} S(a_{(1)})c \ot S(a_{(2)})b
 \qquad\quad\text{and}\qquad\quad
	\sum_{(a)} cS(a_{(1)})\ot  bS(a_{(2)})$$
are well-defined in $A\ot A$. Therefore $(S\ot S)\Delta(a)\in M(A\ot A)$ for all $a$. 
\snl
iii) We would also really like to have $\Delta(S(a))=\sigma(S\ot S)\Delta(a)$ for all $a$, without the need to multiply with $E$. This will follow if we knew that already $(S\ot S)(E)=\sigma(E)$. We will be able to show this in the regular case (see Proposition 4.6 in Section 4), but we have not found an argument in the non-regular case. This seems to be one of the peculiarities in the general situation. 
\hfill$\square$\einspr

We will give more remarks at the end of this section. Let us now consider some applications of the results about the antipode that we obtained.
\snl
In the theory of multiplier Hopf algebras, there is a simple argument to show that the underlying algebra must have local units (cf.\ Proposition 2.2 in [D-VD-Z]). We will now see what happens if we try to generalize this argument here.
\snl
Then the following is obtained.

\inspr{3.9} Proposition \rm
For all $a\in A$ we have that
$$\align
&a\varepsilon_s(A)\subseteq aA \qquad\qquad\qquad \varepsilon_s(A)a\subseteq Aa, \tag"(3.4)"\\
&a\varepsilon_t(A)\subseteq aA \qquad\qquad\qquad \varepsilon_t(A)a\subseteq Aa. \tag"(3.5)"
\endalign$$

\snl\bf Proof\rm:
i) We will first prove the first formula in (3.4). Take $a\in A$ and assume that $\omega$ is a linear functional on $A$ so that $\omega(ab)=0$ for all $b\in A$. For any $b,c,d\in A$, we have that 
$$((S\ot \iota)\Delta(b))(c\ot d)=((S\ot\iota)(\Delta(b)(1\ot d))(c\ot 1)\subseteq S(A)A\ot A\subseteq A\ot A$$
and therefore
$$(\omega\ot\iota)((a\ot 1)((S\ot \iota)\Delta(b))(c\ot d))=0.\tag"(3.6)"$$
Because 
$$(a\ot 1)(S\ot \iota)\Delta(b)\in A\ot A$$
(as this is the same as $R_2(a\ot b)$), we can cancel $d$ in (3.6) and if we write 
$$(a\ot 1)(S\ot \iota)\Delta(b)=\sum_i p_i\ot q_i$$
 with $p_i$ and $q_i$ in $A$, we find that 
$$\sum_i\omega(p_i\,\cdot\,)\ot q_i=0 \qquad\quad \text{in } A'\ot A$$
where $A'$ is the space of all linear functionals on $A$. If we apply the evaluation map, we find that $\sum_i\omega(p_iq_i)=0$. We see that $\sum_{(b)}\omega(aS(b_{(1)})b_{(2)})=0$. As this holds for all $b\in A$, we have shown that $a\varepsilon_s(A)\subseteq aA$. 
\snl
ii) The proof of the first formula in (3.5) is slightly different. Now we start with $a\in A$ and assume again that $\omega(ab)=0$ for all $b\in A$. Again we take $b,d\in A$ and now we allow $c\in M(A)$. Because 
$$\Delta(b)(c\ot d)\subseteq Ac\ot A \subseteq A\ot A,$$
we have now also
$$(\omega\ot\iota)((a\ot 1)(\Delta(b)(c\ot d))=0.\tag"(3.7)"$$
Write $(a\ot 1)\Delta(b)=\sum_i p_i\ot q_i$ with $p_i$ and $q_i$ in $A$. Then we obtain that $\sum_i\omega(p_ic)q_i=0$ for all $c\in M(A)$. Now we apply the antipode and then the evaluation map on $M(A)'\ot A$ (where now $M(A)'$ is the space of all linear functionals on $M(A)$). This gives 
$$\sum_i\omega(p_iS(q_i))=\sum_{(b)}\omega(ab_{(1)}S(b_{(2)}))=0.$$
As this holds for all $b$, we find that $a\varepsilon_t(A)\subseteq aA$.
\snl
iii) The two remaining formulas are obtained in a similar way.
\hfill$\square$\einspr

We know that in the case of a multiplier Hopf algebra, the above result shows that for any element $a$, there exist elements $e$ and $f$ so that $ae=a$ and $fa=a$. And this implies that $A$ has local units (see [Ve] and also [VD-Ve]).
\snl
Unfortunately, we are not able to show that in this general case, the result is still true for weak multiplier Hopf algebras. What we can show is that the result of Proposition 3.9 will guarantee the existence of local units if there are local units for the algebras $\varepsilon_s(A)$ and $\varepsilon_t(A)$. The argument goes as follows.
\snl
Consider two elements $a,a'\in A$ and use that 
$$a'a=\sum_{(a)}a'a_{(1)}S(a_{(2)})a_{(3)}.$$
Because we assume that $\varepsilon_s(A)$ has local units, we find an element $p\in \varepsilon_s(A)$ so that 
$$\sum_{(a)}a'a_{(1)}S(a_{(2)})a_{(3)}=\sum_{(a)}a'a_{(1)}S(a_{(2)})a_{(3)}p$$
and so $a'a=a'ap$. By Proposition 3.9, we get an element $e\in A$ so that $a'a=a'ae$. Because $A^2=A$ we find for all $a$ in $A$ an element $e\in A$ so that $a=ae$. Similarly on the other side and so $A$ has local units.
\snl
In Section 4, where we treat the regular case, we will be able to prove this result and so, we will find that any regular weak multiplier Hopf algebra has local units. 
\snl
In fact, in that section, we will find more results that we have not been able to show in the non-regular case in this section. There is not only the problem of the existence of local units. There is also the problem with the formula for $\Delta(S(a))$ as we mentioned in Remark 3.8.iii. We need to have $(S\ot S)E=\sigma E$, again a formula that we have not been able to show in the non-regular case. Finally, there is the problem of finding the idempotents $F_1$ and $F_2$ as we have indicated already in the beginning of Section 2. We have the feeling that all these properties are related. But we have no idea whether or not they will be true in the non-regular case. See also Section 5 where we suggest some further research here. 
%
%
\nl\nl



\bf 4. Regular weak multiplier Hopf algebras \rm
\nl
In this section, we study {\it regular weak multiplier Hopf algebras}. We recall the definition  we already announced at the end of Section 1 (see Definition 4.1 below). The main result we obtain here is that a weak multiplier Hopf algebra $(A,\Delta)$ is regular if and only if the antipode $S$ maps $A$ into $A$ (and not just into the multiplier algebra $M(A)$) and that it is a bijection. In fact, we only prove one direction in this section (cf.\ Proposition 4.3) while the other direction is formulated (Theorem 4.10) but proven in the appendix. The reason for doing so will be explained. 
\snl
We will also show that some of the peculiarities in the non-regular case completely disappear in the regular case (and therefore, this case is much better understood). In particular, we obtain nicer formulas involving the different idempotents determining the kernels of the canonical maps. As a 'byproduct', we obtain that the underlying algebra will have local units. 
\snl
In this context, we treat the $^*$-case as well.
\snl
The starting point is different from the approach in our first paper [VD-W3]. We use a different (but equivalent) definition.  We introduced this different approach mainly for motivational reasons. The treatment here is mostly self-contained and essentially not based on the results of [VD-W3]. This is completely in accordance with the spirit of this (and the first) paper as we explained it in the introduction already.
\snl
Now recall the definition. Compare with Proposition 4.12 in [VD-W3].

\inspr{4.1} Definition \rm
Let $(A,\Delta)$ be a weak multiplier Hopf algebra (as in Definition 1.14). It is called {\it regular} if the coproduct is regular (as in Definition 1.1) and if also $(A,\Delta^{\text{cop}})$, obtained by flipping the coproduct, satisfies the assumptions of a weak multiplier Hopf algebra.
\hfill$\square$\einspr

Remark that this is equivalent with the requirement that $(A^{\text{op}},\Delta)$ is also a weak multiplier Hopf algebra (where $A^{\text{op}}$ is the algebra $A$ but with the opposite product). Working with  $(A^{\text{op}},\Delta)$ is slightly simpler than with $(A,\Delta^{\text{cop}})$ but it has the disadvantage that we then need to work with a different product on $A$ while the notation for a product is implicit. However we mostly will work in the setting of $(A^{\text{op}},\Delta)$ but avoid the problem by only using formulas written with the original product.
\snl
In what follows, we assume that $(A,\Delta)$ is a regular weak multiplier Hopf algebra. 
\nl
\it Bijectivity of the antipode \rm 
\nl
We will now first of all show that the antipode maps $A$ into $A$ and that it is bijective. We will do this by showing that the antipode $S'$ for the weak multiplier Hopf algebra $(A^{\text{op}},\Delta)$ is the inverse of $S$ (as expected).
\snl
The following result gives the expected equality of the idempotent $E$ we have for the original weak multiplier Hopf algebra $(A,\Delta)$ and the corresponding idempotent for the weak multiplier Hopf algebra $(A^{\text{op}},\Delta)$, obtained by flipping the product. 
\snl
We consider not only the canonical maps $T_1$ and $T_2$, but also the maps $T_3$ and $T_4$, defined from $A\ot A$ to itself by (cf.\ Section 1)
$$T_3(a\ot b)=(1\ot b)\Delta(a)
	\qquad\quad\text{and}\qquad\quad
		T_4(a\ot b)=\Delta(b)(a\ot 1).$$
Remark that these are the maps $T_1$ and $T_2$ for the new pair $(A^{\text{op}} ,\Delta)$. We get the following result for $T_3$ and $T_4$, similarly as what we have for $T_1$ and $T_2$ by definition.

\inspr{4.2} Proposition \rm In the case of a regular weak multiplier Hopf algebra, we also have 
$$(1\ot A)\Delta(A)=(A \ot A)E
\qquad\quad\text{and}\quad\qquad
\Delta(A)(A\ot 1)=E(A\ot A).$$

\snl \bf Proof\rm:
As we assume that the pair $(A^{\text{op}},\Delta)$ also satisfies the conditions of a weak multiplier Hopf algebra (as in Definition 1.14) we have a multiplier $E'$ in $M(A\ot A)$ such that 
$$ (A\ot A)E'=T_3(A\ot A)
\qquad\quad\text{and}\quad\qquad
E'(A\ot A)=T_4(A\ot A).
$$
As we can see from the proof of Proposition 1.6, also now this will imply that $E'$ is the smallest idempotent in $M(A\ot A)$ satisfying 
$$E'\Delta(a)=\Delta(a)
	\qquad\quad\text{and}\qquad\quad
		\Delta(a)E'=\Delta(a) $$
for all $a\in A$. Of course, if we combine this with the result in Proposition 1.6, we find that $E=E'$. This proves the proposition.
\hfill$\square$\einspr

Remark that the canonical idempotent for $(A,\Delta^{\text{cop}})$ is $\sigma E$, obtained from $E$ by applying the flip map. The reader should have this in mind when trying to understand some of the arguments further. The general idea is that formulas we have obtained in the general case (with respect to the maps $T_1$ and $T_2$), have their analogues for the two other canonical maps $T_3$ and $T_4$ in the regular case. Such formulas can be found by either replacing $A$ by $A^{\text{op}}$ or by replacing $\Delta$ by $\Delta^{\text{cop}}$.
\snl
More precisely, if we replace the pair $(A,\Delta)$ by the new pair $(A^{\text{op}},\Delta)$ where $A^{\text{op}}$ is the algebra $A$ but with the opposite product and with the same coproduct, the maps $T_1,T_2$ for the new pair are the maps $T_3,T_4$ for the old one.  We see  that the multiplier $E$ does not change. On the other hand, if we replace the original pair 
$(A,\Delta)$ by $(A,\Delta^{\text{cop}})$, things are slightly more complicated. The  maps $T_1,T_2$ for the new pair are now  respectively $\sigma T_4\sigma, \sigma T_3\sigma$ for the old pair. This implies that $E$ will become $\sigma E$ for the new pair. As before, we use $\sigma$ for the flip map.
\snl
One can also consider the transition from $(A,\Delta)$ to  $(A^{\text{op}},\Delta^{\text{cop}})$. Then the pair $(T_1,T_2)$ is replaced by the pair $(T_2, T_1)$. The idempotent $E$ does not change and also the antipode remains the same in this case.
\nl
Now, we prove the (first) main result for regular weak multiplier Hopf algebras.

\inspr{4.3} Proposition \rm 
Let $(A,\Delta)$ be a regular weak multiplier Hopf algebra. Denote the antipode of $(A,\Delta)$ with $S$ as usual and let $S'$ be the antipode of $(A^{\text{op}},\Delta)$. Then $S$ and $S'$ map $A$ to $A$ and they are each others inverses. 

\snl\bf Proof\rm:
Take $a,b\in A$. Then we have
$$\align 
E(a\ot b)&=\sum_{(a)} a_{(1)}\ot  a_{(2)} S(a_{(3)})b \tag"(4.1)"\\
E(a\ot b)&=\sum_{(b)} b_{(2)} S'(b_{(1)})a \ot  b_{(3)}.\tag"(4.2)"
\endalign$$
The first formula is a rewriting of $T_1R_1(a\ot b)=E(a\ot b)$ while the second one is nothing else but the equation $T_4R_4(a\ot b)=E(a\ot b)$ (which is $T_2R_2(a\ot b)=(a\ot b)E$ but written for $(A^{\text{op}},\Delta)$).  
\snl
If we apply $S$ on the first factor of the right hand sides of the formulas (4.1) and (4.2) and then multiply, we find
$$\align S(a)b
&=\sum_{(a)} S(a_{(1)}) a_{(2)} S(a_{(3)})b\\
&=\sum_{(b)} S(b_{(2)} S'(b_{(1)})a) b_{(3)}\\
&=\sum_{(b)} S(S'(b_{(1)})a)S(b_{(2)}) b_{(3)}
\endalign$$
for all $a,b\in A$. Similarly we have 
$$\align 
(a\ot b)E&=\sum_{(b)} aS(b_{(1)})b_{(2)}\ot b_{(3)} \tag"(4.3)"\\
(a\ot b)E&=\sum_{(a)} a_{(1)}\ot b S'(a_{(3)})a_{(2)}\tag"(4.4)"
\endalign$$
for all $a,b\in A$. Here the first formula is $T_2R_2(a\ot b)=(a\ot b)E$ whereas the second one is $T_3R_3(a\ot b)=(a\ot b)E$ (which is again the formula (4.1), but for the opposite product).
\snl
Now we apply $S$ on the second factor of the right hand sides of the formulas (4.3) and (4.4) and multiply. Then we obtain
$$\align aS(b)
&=\sum_{(b)} aS(b_{(1)}) b_{(2)} S(b_{(3)})\\
&=\sum_{(a)} a_{(1)} S(b S'(a_{(3)})a_{(2)})\\
&=\sum_{(a)} a_{(1)} S(a_{(2)})S(b S'(a_{(3)}))
\endalign$$
for all $a,b\in A$.
\snl
We now combine the two results. Take any $a,b,c\in A$. Applying the first formula for $S(a)b$ and multiplying with $S(c)$, we find
$$S(a)bS(c)=\sum_{(b)} S(S'(b_{(1)})a) S(b_{(2)}) b_{(3)}S(c).\tag"(4.5)"$$
Next we apply the second formula, but with $b_{(3)}$ in the place of $a$ and with $c$ in the place of $b$ and we replace the product $b_{(3)}S(c)$ in the formula (4.5). We get 
$$\align S(a)bS(c)
&=\sum_{(b)} S(S'(b_{(1)})a) S(b_{(2)}) b_{(3)}S(b_{(4)})S(cS'(b_{(5)})) \\
&=\sum_{(b)} S(S'(b_{(1)})a) S(b_{(2)})S(c S'(b_{(3)})) \\
&=\sum_{(b)} S(cS'(b_{(3)})b_{(2)}S'(b_{(1)})a) \\
&=S(cS'(b)a).
\endalign$$
Hence, we obtain $S(a)bS(c)=S(cS'(b)a)$ for all $a,b,c\in A$.
\snl
The left hand side belongs to $A$.  From Proposition 3.6 (applied to $(A^{\text{op}},\Delta)$), we know that any element in $A$ is a linear combination of elements of the form $cS'(b)a$. It follows from this formula that $S(A)\subseteq A$. In fact, we see that $S(A)AS(A)$ is contained in $S(A)$ and again using Proposition 3.6, now applied to $(A,\Delta)$, we see that also $A\subseteq S(A)$ so that $S(A)=A$. Finally, the formula gives $S(a)bS(c)=S(a)S(S'(b))S(c)$. Now we can cancel $S(a)$ and $S(c)$ and find $b=S(S'(b))$ for all $b$. Similarly (or by replacing $(A,\Delta)$ by $(A^{\text{op}},\Delta)$) we find $S'(S(b))=b$ for all $b$. This completes the proof.
\hfill$\square$\einspr

We are using the Sweedler notation in the above proof, but it can be seen from earlier remarks made that this is all justified, i.e.\ that necessary coverings exist.
\snl
As mentioned, we can view the formula (4.4) as (4.1) with $A$ replaced by $A^{\text{op}}$. Similarly (4.2) will then become (4.3). However, we can also interpret the formula (4.2) as the formula (4.1) for $(A,\Delta^{\text{cop}})$, whereas the formula (4.3) turns into the formula (4.4) when replacing the original pair $(A,\Delta)$ by $(A,\Delta^{\text{cop}})$.  
\snl
The result of Proposition 4.3 will also imply that the antipode of the pair $(A,\Delta^{\text{cop}})$ is again $S^{-1}$.
\nl
There are now some important consequences.

\inspr{4.4} Proposition \rm If $(A,\Delta)$ is regular, we have $(S\ot S)E=\sigma E$ where as before $\sigma$ is the flip on $A\ot A$, extended to $M(A\ot A)$. Consequently we then also get
$$\Delta(S(a))=\sigma(S\ot S)\Delta(a)$$
for all $a$.

\snl\bf Proof\rm:
Denote $E'=\sigma(S\ot S)E$.
\snl
Let $a,b$ in $A$. By Proposition 3.7 we have that
$$
\Delta(S(a)b)=\Delta(S(a))\Delta(b)
=(\sigma(S\ot S)\Delta(a))\Delta(b).
$$
Because $S$ is an anti-homomorphism, we get
$$E'\sigma(S\ot S)\Delta(a)=\sigma(S\ot S)(\Delta(a)E)=\sigma(S\ot S)\Delta(a)$$
for all $a$. Hence
$E'\Delta(S(a)b)=\Delta(S(a)b)$
for all $a,b$. Now we have seen in Proposition 3.6 that $A$ is the linear span of elements $S(a)b$ with $a,b\in A$. Therefore we also have
$E'\Delta(a)=\Delta(a)$ for all $a$. Similarly one can show that also $\Delta(a)E'=\Delta(a)$. This implies that $E\leq E'$ by Proposition 1.6.
\snl
Now we use that $(A,\Delta)$ is regular. If we replace $A$ by $A^{\text{op}}$, we know that the antipode is replaced by its inverse. Because  $E$ remains the same, the above result will give $E\leq E''$ where $E''=\sigma(S^{-1}\ot S^{-1})E$. We can apply $\sigma(S\ot S)$ to the inequality $E\leq E''$ because it is an anti-isomorphism and this will give $E'\leq E$. Because we had already that $E\leq E'$, it follows that $E=E'$.
\snl
This proves the first statement. Now, using Proposition 3.7 again and the extra information that $E=\sigma(S\ot S)E$, we can easily obtain that now
$$\Delta(S(a))=\sigma(S\ot S)\Delta(a)$$
for all $a$.
\hfill$\square$\einspr
\nl
\it The idempotents $F_1$, $F_2$, $F_3$ and $F_4$ \rm
\nl
In Section 1, we have proven the existence of idempotent linear maps $G_1$ and $G_2$ and we have assumed in Definition 1.14 that they determine the kernels of the canonical maps $T_1$ and $T_2$. We have a similar result for $(A^{\text{op}},\Delta)$. This gives idempotent maps $G_3$ and $G_4$, giving the kernels of the two other canonical maps $T_3$ and $T_4$. In Remark 1.12, we have seen that it is expected that these maps are given by idempotents in the appropriate multiplier algebras. We now prove this result for regular weak multiplier Hopf algebras in the following proposition. In the formulation, we consider the algebras $A$ and $A^{\text{op}}$, but as mentioned already before, the product we use in the formulas is the original one in $A$.

\inspr{4.5} Proposition \rm
Let $(A,\Delta)$ be a regular weak multiplier Hopf algebra. Then there exists a right multiplier $F_1$ of $A\ot A^{\text{op}}$ and a left multiplier $F_2$ of $A^{\text{op}}\ot A$  so that
$$G_1(a\ot b)=(a\ot 1)F_1(1\ot b)
	\qquad\quad\text{and}\quad\qquad
		G_2(a\ot b)=(a\ot 1)F_2(1\ot b)
$$
for all $a,b\in A$ where $G_1$ and $G_2$ are the maps from $A\ot A$ to $A\ot A$ as defined in Proposition 1.11. Similarly,
there exists a left multiplier $F_3$ of $A\ot A^{\text{op}}$ and a right multiplier $F_4$ of $A^{\text{op}}\ot A$ such that 
$$G_3(a\ot b)=(1\ot b)F_3(a\ot 1)
\quad\quad\text{and}\quad\quad
G_4(a\ot b)=(1\ot b)F_4(a\ot 1)$$
for all $a,b\in A$ where $G_3$ and $G_4$ are the maps $G_1$ and $G_2$ for the pair $(A^{\text{op}},\Delta)$.

\snl\bf Proof\rm:
By definition, we have 
$$(G_1\ot \iota)(\Delta_{13}(a')\Delta_{13}(a)(1\ot b\ot c))= \Delta_{13}(a')\Delta_{13}(a)(1\ot E)(1\ot b\ot c)$$
for all $a,a',b,c\in A$ (see Proposition 1.11). Multiply this equation with an element $d\in A$ from the left in the third factor. By the result obtained in Proposition 4.2, we can replace $(1\ot d)\Delta(a')$ by $(r\ot s)E$ with $r,s\in A$. Using that $E\Delta(a)=\Delta(a)$, we find
$$\align (G_1\ot \iota)((r\ot 1\ot s)\Delta_{13}(a)(1\ot b\ot c))
	&= (r\ot 1\ot s)\Delta_{13}(a)(1\ot E)(1\ot b\ot c)\\
	&=(r\ot 1\ot s)(G_1\ot \iota)(\Delta_{13}(a)(1\ot b\ot c))
\endalign$$
for all $r,s,a,b,c$. We can now cancel $s$ and use that $\Delta$ is full, or also cancel $c$ and apply the counit on the last factor. Then we get $G_1(ra\ot b)=(r\ot 1)G_1(a\ot b)$ for all $r,a,b$. If we combine this with $G_1(a\ot bs)=G_1(a\ot b)(1\ot s)$ for all $a,b,s$ (see a remark following Proposition 1.11), we get the existence of the right multiplier $F_1$ of $A\ot A^{\text{op}}$  satisfying  
$$G_1(a\ot b)=(a\ot 1)F_1(1\ot b)$$
for all $a,b\in A$. This gives the first formula in the proposition.
\snl
Similarly, using that $\Delta(A)(A\ot 1)=E(A\ot A)$ (cf.\ Proposition 4.2), we get the left multiplier $F_2$ of $A^{\text{op}}\ot A$ satisfying 
$$G_2(a\ot b)=(a\ot 1)F_2(1\ot b)$$
for all $a,b\in A$. This gives the second formula.
\snl
The other formulas are obtained in a similar way or by replacing $A$ by $A^{\text{op}}$.
\hfill$\square$\einspr

There is some ambiguity in the choices made above in the formulation of the proposition. Indeed, a right multiplier of $A\ot A^{\text{op}}$ is the same as a left multiplier of $A^{\text{op}}\ot A$. We refer to a remark following Proposition 4.7 below for a motivation of our choices here.
\snl
Remark further that the first part of the proposition is already true provided we have
that the coproduct is regular (as defined in Definition 1.1) and that also 
$$(1\ot A)\Delta(A)=(A\ot A)E
\qquad\quad\text{and}\quad\qquad
\Delta(A)(A\ot 1)=E(A\ot A).$$
We do not need all of the assumptions for regularity of the pair $(A,\Delta)$.
\snl
We also get by assumption iii) of Definition 1.14 that
$$\text{Ker}(T_1)=(A\ot 1)(1-F_1)(1\ot A)
\quad\quad\text{and}\quad\quad
\text{Ker}(T_2)=(A\ot 1)(1-F_2)(1\ot A)$$
and
$$\text{Ker}(T_3)=(1\ot A)(1-F_3)(A\ot 1)
\quad\quad\text{and}\quad\quad
\text{Ker}(T_4)=(1\ot A)(1-F_4)(A\ot 1).$$
\nl 
Now that we have proven that the idempotent maps $G_1$ and $G_2$, as defined in Proposition 1.11, can be expressed using the idempotent multipliers $F_1$ and $F_2$, we can rewrite the defining formulas in Proposition 1.11 in terms of these multipliers. We do this also for the multipliers $F_3$ and $F_4$. This gives the following formulas.

\inspr{4.6} Proposition \rm We have
$$\align E_{13}(F_1\ot 1) &= E_{13}(1\ot E) 
	\quad\quad\text{and}\quad\quad 
		(F_3\ot 1)E_{13}=(1\ot E)E_{13}\tag"(4.6)"\\
	(1\ot F_2)E_{13} &=(E\ot 1)E_{13}
	\quad\quad\text{and}\quad\quad 
	 	E_{13}(1\ot F_4)=E_{13}(E\ot 1).\tag"(4.7)"
\endalign$$

\snl\bf Proof\rm:
If we rewrite the defining formula for $G_1$ in Proposition 1.11 using the relation with the idempotent $F_1$ as obtained in Proposition 4.5, we find
$$(\Delta_{13}(a)(F_1\ot 1)(1\ot b\ot c))= \Delta_{13}(a)(1\ot E)(1\ot b\ot c)$$
for all $a,b,c\in A$. We can replace $\Delta(a)$ by $E$ in this formula and if we also cancel $b$ and $c$, we arrive at the first formula in (4.6). Similarly, we find the first formula in (4.7). 
\snl
The second formula in (4.6) can either be found as the analogue of the first formula in (4.6) for the pair $(A^{\text{op}},\Delta)$ in the place of $(A,\Delta)$ or as the analogue of the first formula in (4.7) for the pair $(A,\Delta^{\text{cop}})$ in the place of $(A,\Delta)$. Similarly for the second formula in (4.7).
\hfill$\square$\einspr

We will now consider some {\it more relations} of the idempotents $E$, $F_1$, $F_2$, $F_3$ and $F_4$, involving the antipode.
\snl
We know that the antipode basically converts the map $T_1$ to $T_2$. The same is also true for the generalized inverses $R_1$ and $R_2$. This should imply that the antipode will convert $F_1$ to $F_2$. More precisely, a more accurate application will give not only that $(S\ot S)E=\sigma E$ (a formula we proved already before in Proposition 4.4), but also that $(S\ot S)F_1=\sigma F_2$. There is however a small problem with the last formula as it requires the extension of $S\ot S$, not only to $M(A\ot A)$ but also to the right multipliers of $A\ot A^{\text{op}}$. There is of course a simple way around this, but from the result in the next proposition, we will see that we actually have (two-sided) multipliers so that the problem does not really occur. 
\snl
Before we formulate and prove the next result, observe that the antipode can also be used to get a relation between the maps $R_1$ and $R_2$ on the one hand and the maps $T_3$ and $T_4$ on the other hand. We have considered this relation already in Propositions 2.11 and 2.12 and in Remark 3.13 of [VD-W3]. Indeed we have 
$$ R_1 (\iota \ot S) = (\iota\ot S) T_3 
	\qquad\text{and}\qquad\quad
		R_2 (S\ot \iota)  = (S\ot \iota) T_4.
$$
Therefore, as we also have observed already in [VD-W3], we expect the following formulas given $F_1$ and $F_2$ in terms of $E$ (as well as similar formulas for $F_3$ and $F_4$). We will first formulate and prove the result. Then we will give more comments.

\inspr{4.7} Proposition \rm
Let $(A,\Delta)$ be a regular weak multiplier Hopf algebra. Then we can express the idempotents $F_1, F_2, F_3$ and $F_4$ all in terms of $E$:
$$\align F_1&=(\iota\ot S)E  \qquad\quad\text{and}\qquad\quad F_3=(\iota\ot S^{-1})E \tag"(4.8)"\\
	F_2&=(S\ot \iota)E \qquad\quad\text{and}\qquad\quad F_4=(S^{-1}\ot \iota)E.\tag"(4.9)"
\endalign$$

\snl\bf Proof\rm:
We know that 
$$\Delta_{13}(a)(F_1\ot 1)=\Delta_{13}(a)(1\ot E)$$
for all $a\in A$. If we apply $S^{-1}$ to the third factor and multiply, we find
$$\sum_{(a)}(S^{-1}(a_{(2)})a_{(1)}\ot 1)F_1=
	\sigma((\iota\ot S^{-1})E)(\sum_{(a)}S^{-1}(a_{(2)})a_{(1)}\ot 1).$$
We know that $\sigma E=(S\ot S)E$ (see Proposition 4.4) and so we have $\sigma((\iota\ot S^{-1})E)=(\iota\ot S)E$. Furthermore, we have 
$$\sum_{(a)}S^{-1}(a_{2})a_{(1)}=S^{-1}(\sum_{(a)}S(a_{1})a_{(2)})$$
for all $a$. As 
$$(1\ot a)E=\sum_{(a)}S(a_{(1)})a_{(2)}\ot a_{(3)},$$
we see that $\sum_{(a)}S(a_{(1)})a_{(2)}$ belongs to the first leg of $E$ and so that $\sum_{(a)}S^{-1}(a_{(2)})a_{(1)}$ belongs to the second leg of $E$ (as the antipode flips the legs of $E$). Now, because the two legs of $E$ commute (as $1\ot E$ commutes with $E\ot 1$), we conclude that  
$$\sum_{(a)}(S^{-1}(a_{(2)})a_{(1)}\ot 1)F_1=
	(\sum_{(a)}S^{-1}(a_{(2)})a_{(1)}\ot 1)((\iota\ot S)E)$$
for all $a$. If we apply this to the second leg of $\Delta(a)$ in the place of $a$ and multiply, we find $(a\ot 1)F_1 =(a\ot 1)((\iota\ot S)E)$ for all $a$ as 
$$\sum_{(a)}a_{(3)}S^{-1}(a_{(2)})a_{(1)}=a$$
for all $a$. Then the  formula for $F_1$ is proven.
\snl
To prove the one for $F_2$, we start from the equation
$$(1\ot F_2)\Delta_{13}(a)=(E\ot 1)\Delta_{13}(a)$$
and use a similar argument.
\snl
The two other formulas are obtained from the previous ones by the standard procedure, either replacing $(A,\Delta)$ by $(A^{\text{op}},\Delta)$ or by $(A,\Delta^{\text{cop}})$.
\hfill$\square$\einspr

We really should have been more careful with the above argument. We could e.g.\ have multiplied the starting formula with an element $1\ot b\ot c$ from the right where $b,c\in A$. Then we get elements in $A\ot A\ot A$ on both sides and we can safely continue. 
\snl
This would also yield a correct interpretation of the formula $F_1=(\iota\ot S)E$ as
$$(a\ot 1)F_1(1\ot S(b))=(\iota\ot S)((a\ot b)E)$$
for all $a,b\in A$. Another possible interpretation is by extending the map $\iota\ot S$ to $M(A\ot A)$ first, but that is essentially the same story. Similarly for the other formulas.
\snl
It now should  also be clear  why we have considered e.g.\ $F_1$ as a right multiplier of $A\ot A^{\text{op}}$ and not as a left multiplier of $A^{\text{op}}\ot A$ (cf.\ the first formula in (4.8)). However, we see from the following corollary of this proposition that it does not really matter.

\inspr{4.8} Proposition \rm
The elements $F_1$, $F_2$, $F_3$ and $F_4$ are all elements in $M(A\ot A^{\text{op}})$ (which is the same as $M(A^{\text{op}}\ot A$). Furthermore we see that not only $(S\ot S)E=\sigma E$ but also $(S\ot S)F_1=\sigma F_2$ and $(S\ot S)F_3=\sigma F_4$. And consequently we get that all these multipliers are left invariant by $S^2\ot S^2$. 
\hfill$\square$\einspr

Before we continue, we will apply our results and prove the following.

\inspr{4.9} Proposition \rm
If $(A,\Delta)$ is a regular weak multiplier Hopf algebra, then $A$ has local units.

\bf \snl  Proof\rm:
Let $a\in A$ and suppose that $\omega$ is a linear functional on $A$ so that $\omega(ab)=0$ for all $b\in A$. We will show that then $\omega(a)=0$. This gives $a\in aA$. In a similar way (or by applying the antipode), one can show that also $a\in Aa$. Then $A$ has local units (see [Ve] or [VD-Ve]).
\snl
So assume that $\omega$ is $0$ on all of $aA$. We have seen in Proposition 3.9 that then also $\omega(ay)=0$ for all $y\in \varepsilon_s(A)$. In Lemma 3.2 we have seen that $\varepsilon_s(A)$ is the left leg of $E$ and in Lemma 3.4 that it is a subalgebra of $M(A)$. Then it follows that
$$(\omega\ot\iota)(a\ot 1)E(y\ot 1))=0$$
and by applying $S$ and using that $F_1=(\iota\ot S)E$ that also
$$(\omega\ot\iota)(a\ot 1)F_1(y\ot 1))=0$$
for all $y\in \varepsilon_s(A)$. Now we have
$$(a\ot 1)F_1(1\ot b)=R_1T_1(a\ot b)=\sum_{(a)}a_{(1)}\ot S(a_{(2)})a_{(3)}b\tag"(4.10)"$$
for all $b$ and it will follow that
$$\sum_{(a)}\omega(a_{(1)}y)S(a_{(2)})a_{(3)}=0 \tag"(4.11)" $$
for all $y\in \varepsilon_s(A)$. However, for any two elements $p,q$ in $A$, we find that 
$$(pq\ot 1)F_1=(p\ot 1)(q\ot 1)F_1=\sum_{(q)}pq_{(1)}\ot S(q_{(2)})q_{(3)}$$
by using the formula (4.10) above
and this belongs to $A\ot \varepsilon_s(A)$. Because $A^2=A$ we will also have that 
$$\sum_{(a)}a_{(1)}\ot S(a_{(2)})a_{(3)}\in A\ot \varepsilon_s(A).$$
Therefore, we can apply the evaluation map (as in the proof of Proposition 3.9) for the equation (4.11) and we get
$$\sum_{(a)}\omega(a_{(1)}S(a_{(2)})a_{(3)})=0.$$
This means $\omega(a)=0$ and the proof is complete.
\hfill$\square$\einspr 

It remains an open problem whether or not the underlying algebra of any weak multiplier Hopf algebra has local units. As we mentioned already at the end of the previous section, this is true as soon as the algebras $\varepsilon_s(A)$ and $\varepsilon_t(A)$ have local units. Indeed, it is a typical problem for weak multiplier Hopf algebras because we know that the result is true for multiplier Hopf algebras, regular or not.
\nl
Before we consider some special cases and examples, we formulate what is to be considered as the {\it main result} of this section. 

\iinspr{4.10} Theorem \rm
Let $(A,\Delta)$ be a weak multiplier Hopf algebra. Then it is regular if and only if the antipode maps $A$ to $A$ and is bijective.
\hfill$\square$\einspr 

One direction has been obtained in Proposition 4.3. The proof of the other direction is given in the appendix. The reason we do not give the proof here, but rather in an appendix, is that the proof seems to be more involved than expected. It needs several properties that are not really related with the topic of this section. 
\snl
The reader can also have a look at the treatment of the regular case in [VD-W3] as there, regular weak multiplier Hopf algebras are defined by the requirement that the antipode is bijective from $A$ to itself. However, we should also have in mind that the definition of a weak Hopf algebra in this first paper is slightly more restrictive than the one we use in this paper, making thing more complicated.
\nl
\it Special cases and examples \rm
\nl
We will not consider the groupoid examples in this place. We have already mentioned that they give {\it regular} weak multiplier Hopf algebras. Now, in view of Theorem 4.10, this is also obvious as the antipode is involutive and so bijective.
\snl
Instead, we will focus  on the case of a weak Hopf algebra. But first, let us now briefly consider the involutive case. Remember that a weak multiplier Hopf algebra $(A,\Delta)$ where $A$ is a $^*$-algebra and $\Delta$ a $^*$-homomorphism is called a weak multiplier Hopf $^*$-algebra (cf.\ a remark following Definition 1.14 in Section 1).
The next result is expected.

\iinspr{4.11} Proposition \rm
If $(A,\Delta)$ is a {\it weak multiplier Hopf $^*$-algebra}, then it is regular. The antipode satisfies $S(S(a)^*)^*=a$ for all $a\in A$. And not only do we have $E^*=E$ but also 
$$F_1^*=F_3 \qquad\quad\text{and}\qquad\quad F_2^*=F_4.$$

\snl\bf Proof\rm:
We know that the coproduct is regular and that 
$$T_3(a^*\ot b^*)=T_1(a\ot b)^*
     \qquad\quad\text{and}\qquad\quad
       T_4(a^*\ot b^*)=T_2(a\ot b)^*. 
$$
From this it easily follows that $(A^{\text{op}},\Delta)$ is again a weak multiplier Hopf algebra. In particular, $(A,\Delta)$ is a regular weak multiplier Hopf algebra. 
\snl
As we know already that the idempotent $E$ is the same for both $(A,\Delta)$ and $(A^{\text{op}},\Delta)$, we conclude from this that $E^*=E$. In fact, this also follows from the fact that $E$ is the smallest idempotent in $M(A\ot A)$ such that $\Delta(a)=\Delta(a)E=E\Delta(a)$ for all $a$ (by taking adjoints).
\snl
As we have seen already before, we will have that the antipode of $(A^{\text{op}},\Delta)$ is $S^{-1}$. But it is also obviously given by $a\mapsto S(a^*)^*$.  This will imply the property of $S$ as in the formulation of the proposition. 
\snl
Finally, the equalities $F_1^*=F_3$ and  $F_2^*=F_4$ follow from the formulas in Proposition 4.7.

\hfill$\square$\einspr 

Next, consider again the case of a weak Hopf algebra. In Section 2, we showed already that any weak Hopf algebra is a weak multiplier Hopf algebra (cf.\ Proposition 2.10). We now prove the following converse result. We refer to the papers [B-N-S] and [N-V2] for some of the notions below.

\iinspr{4.12} Proposition \rm
Let $(A,\Delta)$ be a regular weak multiplier Hopf algebra and assume that $A$ has an identity. Then it is a weak  Hopf algebra. 

\snl\bf Proof\rm:
As the algebra is unital, we have $M(A)=A$ and $M(A\ot A)=A\ot A$. Therefore the coproduct $\Delta$ is a coproduct on $A$ in the usual sense. There is a counit by our assumptions. 
\snl
And of course we have $E=\Delta(1)$ in this case. The behavior of $\Delta$ on $\Delta(1)$ is then one of our assumptions in Definition 1.14.
\snl
The antipode is a map from $A$ to $A$ in this case as $M(A)=A$. We know that $E(a\ot 1)=\sum_{(a)}a_{(1)}\ot a_{(2)}S(a_{(3)})$ and if we apply $(\varepsilon\ot \iota)$ we find the formula
$$(\varepsilon\ot \iota)(\Delta(1)(a\ot 1))=\sum_{(a)}a_{(1)}S(a_{(2)})$$
for all $a\in A$. Similarly we have $(1\ot a)E=\sum_{(a)}S(a_{(1)})a_{(2)}\ot a_{(3)}$ and if we apply $(\iota\ot\varepsilon)$, we arrive at
$$(\iota\ot\varepsilon)((1\ot a)\Delta(1))=\sum_{(a)}S(a_{(1)})a_{(2)}$$
for all $a$. Because we also have $\sum_{(a)}S(a_{(1)})a_{(2)}S(a_{(3)})=S(a)$ for all $a$ we have all the necessary statements about the antipode.
\snl
It remains to prove the extra properties of the counit. They are called {\it weak multiplicativity} of the counit in the original paper [B-N-S] and say
$$\align  \varepsilon(abc)&=\sum_{(b)} \varepsilon(ab_{(2)})\varepsilon(b_{(1)}c)\tag"(4.12)"\\
\varepsilon(abc)&=\sum_{(b)} \varepsilon(ab_{(1)})\varepsilon(b_{(2)}c) \tag"(4.13)"
\endalign$$
for all $a,b,c\in A$. 
To do this, we start with 
$$(1\ot a)\Delta(b)(c\ot 1)=\sum_{(c)}(1\ot a)\Delta(bc_{(1)}(1\ot S(c_{(2)})),$$ 
true for all $a,b,c\in A$. If we apply $\varepsilon\ot\varepsilon$ we find
$$(\varepsilon\ot\varepsilon)((1\ot a)\Delta(b)(c\ot 1))=\sum_{(c)}\varepsilon(abc_{(1)}S(c_{(2)})).$$  
If we apply this with $a=1$ we find 
$$\varepsilon(bc)=(\varepsilon\ot\varepsilon)\Delta(b)(c\ot 1)=\sum_{(c)}\varepsilon(bc_{(1)}S(c_{(2)}))$$
and if we use this formula with $b$ replaced by $ab$ in the previous formula, we find
$$(\varepsilon\ot\varepsilon)((1\ot a)\Delta(b)(c\ot 1))=\varepsilon(abc)$$
for all $a,b,c\in A$. This gives one of the properties of the counit we need, namely the first one (4.12).
\snl
Now we use the assumption that $A$ is regular. Then we can apply the previous result for $(A,\Delta^{\text{cop}})$ and since the counit is the same, we find the other formula (4.13)
$$(\varepsilon\ot\varepsilon)((a\ot 1)\Delta(b)(1\ot c))=\varepsilon(abc)$$
for all $a,b,c\in A$. We also find this formula from the other one, applied to  $(A^{\text{op}},\Delta)$.

\hfill$\square$\einspr

There is something peculiar going on here. We will explain this in a remark below. First, let us consider the finite-dimensional case.

\iinspr{4.13} Proposition \rm
Let $(A,\Delta)$ be a {\it finite-dimensional weak Hopf algebra}. Then it is a regular weak multiplier Hopf algebra. Conversely, if $(A,\Delta)$ is a regular weak multiplier Hopf algebra and if $A$ is finite-dimensional, then it is a finite-dimensional weak Hopf algebra.

\snl\bf Proof\rm:
If $(A,\Delta)$ is a finite-dimensional weak Hopf algebra, it is a weak multiplier Hopf algebra as we have proven in Proposition 2.10. Because for a finite-dimensional weak Hopf algebra, the antipode is bijective (see e.g.\ Theorem 2.10 in [B-N-S]), it follows that we have a regular weak multiplier Hopf algebra. This proves one implication.
\snl
The converse will follow from the previous result. Indeed, because the weak multiplier Hopf algebra is assumed to be regular, it has local units and as the algebra is finite-dimensional, it has to have a unit. Then, as a consequence of the previous proposition, it will be a weak Hopf algebra.
\hfill$\square$\einspr

We finish this section by the following remarks.

\iinspr{4.14} Remark \rm
i) In Proposition 2.10 we have shown that any weak Hopf algebra is a weak multiplier Hopf algebra. A remarkable fact was that we did not use the weak multiplicativity axioms of the counit (as formulated in the formulas (4.12) and (4.13).
\snl
ii) On the other hand, we see that for any weak multiplier Hopf algebra with a unital algebra, we have one part of these axioms, namely formula (4.12). This implies that for any weak Hopf algebra, this axiom will follow from the other ones (not including any of the weak multiplicativity formulas of course). In fact, this is easy to show using arguments of our Proposition 4.12.
\snl
iii) In order to get also the second formula (4.13), we need regularity of the given weak multiplier Hopf algebra. 
\snl
iv) On the other hand, in the finite-dimensional case, regularity follows from the axioms for a weak Hopf algebra.
\hfill$\square$\einspr

As a consequence of all this, one may wonder first of all if our notion is more general  in the finite-dimensional case. Of course, this would only be the case if there are non-regular finite-dimensional weak multiplier Hopf algebras. This is open.
\snl
In the infinite dimensional case, a similar question can be asked. If the {\it algebra is unital}, we have the two following implications:
\snl
i) If $(A,\Delta)$ is a regular weak multiplier Hopf algebra then it is a weak Hopf algebra. In particular, the counit satisfies the two weak multiplicativity axioms (4.12) and (4.13).\newline
ii) If $(A,\Delta)$ satisfies all the axioms of a weak Hopf algebra, except for the weak multiplicativity axiom (4.13) - recall that the other is automatic - then it is a weak multiplier Hopf algebra.
\snl
This seems to suggest that one can remove the weak multiplicativity axioms for the counit. On the other hand, this contradicts the desire to have a set of axioms that is self-dual. Indeed, the formulas (4.12) and (4.13) are dual to the formulas
$$\align (\Delta\ot \iota)E&=(1\ot E)(E\ot 1)\tag "(4.14)"\\
(\Delta\ot \iota)E&=(E\ot 1)(1\ot E).\tag"(4.15)"\endalign$$
\snl
And to make the mystery complete, we refer to the first paper [VD-W3] where we also found that (4.14) was more natural than (4.15). This condition alone is in fact self-dual already. The second one was motivated there also by regularity. And it was taken as an extra assumption for this reason. See Proposition 3.9 and Assumption 3.10 in [VD-W3].
\snl
We will come back to this in Section 5 where we discuss further research.
\newpage


\bf 5. Conclusions and further research \rm
\nl
In this paper, we have given a precise definition of a weak multiplier Hopf algebra $(A,\Delta)$ (cf.\ Definition 1.14 in Section 1). The definition is close in spirit to the one of a multiplier Hopf algebra (as introduced in [VD1]). Indeed, apart from some natural conditions on the coproduct, conditions are formulated in terms of the ranges and the kernels of the canonical maps, given on $A\ot A$ by
$$T_1(a\ot b)=\Delta(a)(1\ot b)
\qquad\quad\text{and}\qquad\quad
T_2(a\ot b)=(a\ot 1)\Delta(b).$$
Recall that these maps are assumed to be bijective in the case of a multiplier Hopf algebra but that this is {\it no longer the case} for weak multiplier Hopf algebras.
\snl
Where the conditions come from has been explained in a preliminary paper [VD-W3]. There we have shown that indeed, from very natural considerations, we arrive at the conditions that we use in the definition of a weak  multiplier Hopf algebra in this paper.
\snl
The main result in this paper is the construction of the antipode and its properties. We gave an alternative definition of a weak multiplier Hopf algebra in terms of the antipode. In the regular case, that is when the antipode is a bijective map from $A$ to itself, we found nicer and more formulas. It is shown that any  weak Hopf algebra (as introduced in [B-N-S]) is a weak multiplier Hopf algebra and that moreover the finite-dimensional weak Hopf algebras are precisely the regular weak multiplier Hopf algebras with a finite-dimensional underlying algebra.
\snl
In [VD-W4] we continue this work and study in greater detail the source and target maps and source and target algebras as up to some extend already introduced in this paper in Section 3. We also give non-trivial examples of (regular) weak multiplier Hopf algebras in [VD-W4]. Related with this topic is the study of the notion of separability for (non-degenerate) algebras without identity in [VD4]. In [VD-W5] we study integrals on weak multiplier Hopf algebras. We include an existence proof of integrals in the finite-dimensional case. This is also the natural setting to consider duality. Furthermore, we are working on a theory of cointegrals for weak multiplier Hopf algebras. In that case, it seems to be possible to show that integrals automatically exist. Finally, the work on weak multiplier Hopf algebras will also be formulated in an 'algebroid' framework later. This will be heavily based on the results on the source and target maps as obtained in [VD-W4].  All of this is work 'in progress'.
\snl
We have shown that any groupoid $G$ gives rise, in two ways, to a weak multiplier Hopf algebra. On the one hand, there is the algebra $K(G)$ of complex functions with finite support in $G$ and pointwise product. The coproduct is dual to the product in $G$. On the other hand, there is the convolution algebra $\Bbb C G$, with the 'pointwise' coproduct. These two cases are in duality with each other. If $G$ is finite, these examples are typical examples of weak Hopf algebras as the underlying algebras have an identity. If however $G$ is no longer assumed to be finite, the function algebra never has an identity while also the convolution algebra will not have an identity if the set of units in $G$ is not finite. Therefore these two cases can not be treated within the framework of weak Hopf algebras and weak multiplier Hopf algebras are necessary.
\snl
There are many other interesting examples of weak multiplier Hopf algebras that are not weak Hopf algebras. By lack of space, we have not included these examples in this paper, but we have done this in our forthcoming papers on the subjects as mentioned already.
\snl
Finally, this work should also be compared with the analytical theory as developed by 'the French school' (M.\ Enock, F.\ Lesieur, M.\ Vallin, ...), see e.g.\ [V] and references therein. In particular, it should be investigated if it is possible, given a weak multiplier Hopf $^*$-algebra with positive integrals, to construct an operator algebra representation that fits into the analytical theory. It is expected that this procedure will work, just as in the case of multiplier Hopf $^*$-algebras. Moreover, as this was also true in that situation, one might expect to gain a better understanding of the analytical theory from our work on weak multiplier Hopf algebras.
\snl
This suggests further research in one direction. In some sense, the further investigation of the non-regular case, is the 'opposite' direction. 
\snl
First of all, one needs to look for examples of weak multiplier Hopf algebras that are not regular. A Hopf algebra with a non-invertible antipode will of course be such an example. But one should consider cases where e.g.\ even the coproduct is not regular. Also, one needs to find examples where the antipode $S$ does not map $A$ to itself. Such examples can not be found within the setting of weak Hopf algebras as there the algebra is unital. We also need to mention that finding such examples is even open for multiplier Hopf algebras. It is not immediately clear if this case is easier or more complicated than the case of weak multiplier Hopf algebras. 
\snl
A next problem to investigate is the existence of local units. We have seen that this is rather a problem about the source and target algebras (see a remark at the end of Section 3).
\snl
There is the problem with the anti-coalgebra property of the antipode (see Proposition 3.7) and with the related equality $(S\ot S)E=\sigma E$.
\snl
Finally, we refer to the discussion at the end of Section 4 about the weak multiplicativity axioms for the counit and the related commutativity of the legs of $E$. It seems worthwhile to investigate if these axioms are really necessary.
%
\nl\nl



\bf Appendix A. Regularity and bijectivity of the antipode. \rm
\nl
In Section 4 (Definition 4.1), we have defined {\it regularity} of a weak multiplier Hopf algebra $(A,\Delta)$ by the requirement that the pair $(A,\Delta^{\text{cop}})$, obtained by flipping the coproduct on $A$, again satisfies the axioms of a weak multiplier Hopf algebra as formulated in Definition 1.14 of this paper. We have remarked that this is equivalent with the requirement that the pair $(A^{\text{op}},\Delta)$, obtained by taking $A$ with the opposite product, but with the original coproduct, also satisfies the conditions of Definition 1.14.
\snl
We have shown that for a regular weak multiplier Hopf algebra, apart from some other nice properties, the antipode $S$ maps $A$ into $A$ (and not just in $M(A)$ as in the general case) and that it is bijective. This is proven in Proposition 4.3. 
\snl
In our first paper on the subject [VD-W3], where we spent a lot of effort to motivate the definition of a weak multiplier Hopf algebra, we not only considered a slightly different definition from what we do in this paper (compare Definition 4.1 of [VD-W3] with Definition 1.14 in this paper), but we also had a different initial notion of regularity. Indeed, in [VD-W3] we have called, for motivational reasons, a weak multiplier Hopf algebra $(A,\Delta)$ regular if the antipode is bijective from $A$ to itself (see Definition 4.6 in [VD-W3]). 
\snl
From Proposition 4.3 in this paper, we know that the regularity as defined in this paper, will give regularity as in the first paper. Also the converse is true, as we announced in Section 4 (see Proposition 4.10). However, we have not yet given a proof there because, as we will see in this appendix, where we do prove this result, it is not so easy to obtain this property. We feel that this is somewhat strange as the result is very easy to show in the case of multiplier Hopf algebras. On the other hand, we have not found a simpler argument than the one we present here, in this appendix.
\nl
So, in what follows, we have a weak multiplier Hopf algebra $(A,\Delta)$ and we assume that the antipode $S$, as obtained in Section 2, maps $A$ into $A$ and  that it is  bijective (that is one to one and onto). Some of the intermediate results are true in more general situations and we will mention this when it happens.

\inspr{A.1} Proposition \rm The coproduct $\Delta$ on $A$ is regular (as in Definition 1.1). When considered on $A^{\text{op}}$, it is still a full coproduct (as in Definition 1.4) with a counit (Definition 1.3).

\snl\bf Proof\rm:
We know that the maps $R_1$ and $R_2$, as introduced in Proposition 2.3, map $A\ot A$ into itself. Using now that the antipode is an anti-homomorphism from $A$ to itself, we easily find from the formulas in Proposition 2.4 that
$$R_1(a\ot S(b))=(\iota \ot S)((1\ot b)\Delta(a))
\quad\quad\text{and}\quad\quad
R_2(S(a)\ot b)=(S\ot\iota)(\Delta(b)(a\ot 1))$$
for all $a,b\in A$. Because $S$ is bijective, if follows that
$$(1\ot b)\Delta(a)
\qquad\qquad\text{and}\qquad\qquad
\Delta(b)(a\ot 1)$$
belong to $A\ot A$ for all $a,b\in A$. This precisely means that the coproduct is regular.
\snl
We have mentioned already in Section 1 (see a remark following Definition 1.1) that coassociativity, as formulated in Definition 1.1, remains true on the opposite algebra.
\snl
That the coproduct $\Delta$ is still full on the opposite algebra $A^{\text{op}}$ has been argued in Proposition 1.7 of [VD-W1]. See also a remark following Lemma 1.11 in [VD-W3]. 
\snl
Finally it is easy to see that the original counit is also the counit for the coproduct on the opposite algebra. Again see a remark after Definition 1.3.
\hfill$\square$\einspr

For the following result, it is in fact sufficient that the coproduct is regular. Recall the maps $\varepsilon_s$ and $\varepsilon_t$ as given in Definition 3.1 and Lemma 3.2 where it is shown that their images are (in a certain sense) the left and right legs of $E$ respectively. 

\inspr{A.2} Proposition \rm
For all $a\in A$ we have
$$E(a\ot 1)\in A\ot\varepsilon_t(A)
\qquad\quad\text{and}\qquad\quad
(1\ot a)E\in \varepsilon_s(A)\ot A.$$

\snl\bf Proof\rm:
We know that 
$$E(a\ot b)=T_1R_1(a\ot b)=\sum_{(a)}a_{(1)}\ot a_{(2)}S(a_{(3)})b$$
for all $a,b\in A$. Then we see that 
$$E(aa'\ot 1)=\sum_{(a)}a_{(1)}a'\ot \varepsilon_t(a_{(2)})$$
for all $a,a'\in A$. The right hand side belongs to $A\ot \varepsilon_t(A)$ because the coproduct is regular so that
$$\sum_{(a)}a_{(1)}a'\ot a_{(2)}$$
is in $A\ot A$ for all $a,a'\in A$. Finally we use that $A^2=A$ implying that also 
$$E(a\ot 1)\in A\ot\varepsilon_t(A)$$
for all $a\in A$. 
\snl
The second formula is proven in a completely similar way.
\hfill$\square$\einspr

The formulas in the next lemma will need to be explained. This is done in the proof and commented in a remark following the result. In particular, we use $m$ for the multiplication map from $A\ot A$ to $A$, but we will need some extension of this map to a larger space.

\inspr{A.3} Lemma \rm
$$m(S\ot \iota)E=1
\qquad\quad\text{and}\qquad\quad
m(\iota \ot S)E=1.$$

\snl\bf Proof\rm:
When $a,b\in A$, we have
$$E(a\ot b)=\sum_{(a)}a_{(1)}\ot a_{(2)}S(a_{(3)})b$$
in $A\ot A$ and we can apply $m(S\ot\iota)$. On the right hand side, we get
$$\sum_{(a)}S(a_{(1)})a_{(2)}S(a_{(3)})b=S(a)b$$
by Proposition 2.9. When we apply it on the left hand side, we can consider this as 
$$S(a)(m(S\ot\iota)E)b.$$
This means that we extend multiplication in a certain way and we arrive at (an interpretation of) the first formula in the lemma.
\snl
In a completely similar way, we get and interpret the second formula.
\hfill$\square$\einspr

First remark that we do not need the extra assumptions on the antipode. However, if the coproduct is regular (as is true when we have our extra assumptions), we have seen that $E(a\ot 1)\in A\ot\varepsilon_t(A)$ and because $\varepsilon_t(A)$ is in $M(A)$, we can already apply $m(S\ot\iota)$ to this equation to get
$$S(a)(m(S\ot \iota)E)=S(a)$$
for all $a$. This is now an equation in $M(A)$ (so if regularity of the coproduct is assumed). In some sense, this result is a bit stronger, although this is only a valid statement as far as the interpretation is concerned. We will see  that the weaker interpretation is sufficient for the application of this result in the next proposition (and further). 

\inspr{A.4} Proposition \rm
For all $a\in A$ we  have 
$$S(\varepsilon_t(a))=\varepsilon_s(S(a))
\qquad\quad\text{and}\qquad\quad
S(\varepsilon_s(a))=\varepsilon_t(S(a)).$$

\snl\bf Proof\rm:
i) First observe that the source and target maps have their images in $M(A)$ and that we can extend the antipode to the multiplier algebra (see a remark following Proposition 3.6). Therefore the left hand sides are defined in $M(A)$. This is also true for the right hand sides because we assume that $S(A)\subseteq A$.
\snl
Next we take $a\in A$  and we start with the formula
$$\Delta(S(a))=(\sigma(S\ot S)\Delta(a))E,$$
proven in Proposition 3.7. Multiply this formula with $p\ot q$ from the left, where $p,q$ are elements $A$. We arrive at a formula in $A\ot A$ and we can apply $m(\iota\ot S)$. 
\snl
The left hand side will be $p\varepsilon_t(S(a))S(q)$. For the right hand side, we can use the formula
$$m(\iota\ot S)((r\ot s)E)=rS(s)$$
for $r,s\in A$ as obtained in the previous lemma. We get
$$\sum_{(a)}pS(a_{(2)})S(S(a_{(1)}))S(q)=pS(\varepsilon_s(a))S(q).$$
As this is true for all $p,q\in A$, and because $S(A)=A$, we can cancel $p$ and $S(q)$ and we get the first formula.
\snl
The second formula is proven in a completely similar way.
\hfill$\square$\einspr

Remark that we only need that $S$ maps $A$ into $A$ and nothing more in order to formulate and prove the above result. Indeed, for the last step in the proof, we can also use that and $S(A)A=A$ (see Proposition 3.6), so that we still can cancel $p$ and $S(q)$. In the more general setting, it would be harder to interpret the formulas (although the argument will still be useful).
\snl
The following lemma is again true without any further assumption on the antipode or the coproduct.

\inspr{A.5} Lemma \rm
For $y\in \varepsilon_s(A)$ and $x\in \varepsilon_t(A)$ we have
$$E(y\ot 1)=E(1\ot S(y))
\qquad\quad\text{and}\qquad\quad
(1\ot x)E=(S(x)\ot 1)E.$$

\bf\snl Proof\rm: 
Take $y\in \varepsilon_s(A)$ and $a\in A$. From Lemma 3.3 we know that \newline 
$\Delta(ya)=(1\ot y)\Delta(a)$. Then it follows that
$$\align E(ya\ot 1)
&=\sum_{(ya)}(ya)_{(1)}\ot (ya)_{(2)}S((ya)_{(3)})\\
&=\sum_{(a)}a_{(1)}\ot a_{(2)}S(ya_{(3)})\\
&=\sum_{(a)}a_{(1)}\ot a_{(2)}S(a_{(3)})S(y)\\
&=E(a\ot S(y)).
\endalign$$
Because this holds for all $a$, we find $E(y\ot 1)=E(1\ot S(y))$.
\snl
The proof of the other formula is similar.
\hfill$\square$\einspr

Remark that this result is true for all $y\in A_s$ and all $x\in A_t$ where $A_s$ and $A_t$ are the source and target algebras respectively, as studied further in [VD-W4]. However, we only need the result for the images in these algebras of the source and target maps as above.
\nl
We recall the following notation (introduced already in the proof of Proposition 4.4). 

\inspr{A.6} Notation \rm
We let $E'=\sigma(S\ot S)E$.
\hfill$\square$\einspr

We have explained that $S\ot S$ can be extended to the multiplier algebra $M(A\ot A)$ and of course the same is true for the flip map $\sigma$. Therefore $E'$ is well-defined in $M(A\ot A)$. It also follows that $E'$ is again an idempotent. 
\snl
We will show that actually $E'=E$, as expected because we think of $E$ as $\Delta(1)$. In fact, this seems to be the main point in order to get the result we want to prove in this appendix. And it is somewhat surprising that the proof of this is rather involved. 
\snl
Before we can show this, we first draw some immediate conclusions from the two previous results. 

\inspr{A.7} Lemma \rm
For $y\in \varepsilon_s(A)$ and $x\in \varepsilon_t(A)$ we also have
$$E'(y\ot 1)=E'(1\ot S(y))
\qquad\quad\text{and}\qquad\quad
(1\ot x)E'=(S(x)\ot 1)E'.$$

\vskip -0.5cm
\hfill$\square$\einspr

This is proven by applying the map $\sigma(S\ot S)$ to the equations in Lemma A.5 and using the result of Proposition A.4.
\snl
Now we are ready to prove that $E'=E$, in other words that $\sigma(S\ot S)E=E$. We now need again that $S$ is bijective from $A$ to itself.

\inspr{A.8} Proposition \rm
We have $E'=E$.

\snl\bf Proof\rm:
i) As in the proof of Proposition 4.4, we find that 
$E\leq E'$ in the sense that $E'E=EE'=E$.
\snl
ii) We will now show that also $EE'=E'$. This will complete the proof.
\snl
Because the right leg of $E$ is $\varepsilon_t(A)$ in the sense of Lemma 3.2, an application of the formula $(1\ot x)E'=(S(x)\ot 1)E'$ for $x\in \varepsilon_t(A)$, proven in Lemma A.7, will give 
$$E_{13}(1\ot E')=(F_1\ot 1)(1\ot E')$$
where now $F_1$ is used to denote $(\iota\ot S)E$.
We can apply multiplication on the first two factors and use that $mF_1=m(\iota\ot S)E=1$ (as shown in Lemma A.3) to obtain that $EE'=E'$. The argument is made precise by multiplying  with  $u\ot 1\ot v$ from the left and with $1\ot p\ot q$ from the right, where of course $u,v,p,q\in A$. Doing this, we only need the weaker interpretation of the formula $m(\iota\ot S)E=1$ as explained in the proof of Lemma A.3.
\hfill$\square$\einspr

Before we continue, we want to include some more remarks.

\inspr{A.9} Remarks \rm
i) The formula $E\leq E'$ is true in general, so without the extra assumption on $S$. 
\snl
ii) Because $E\leq E'$, we see that Lemma A.7 generalizes Lemma A.5.
\snl
iii) If we start with the assumption that $(A,\Delta)$ is regular, we find this formula also for $A^{\text{op}}$ and as the antipode is then replaced by its inverse, this formula is essentially transformed to $E'\leq E$. So we get $E'=E$. This is what we used in the proof of Proposition 4.4.
\snl
iv) In the converse setting, we see that it is still true that $E'=E$, but it turns out to need more work to show this.
\hfill$\square$\einspr

We will now continue the proof of our result, essentially by reversing the arguments in the treatment of the direct result in Section 4. We {\it define} the elements $F_1$ and $F_2$ by the formulas in Proposition 4.7 and we show that they satisfy the properties in Proposition 4.6 and that they give the maps $G_1$ and $G_2$ as in Proposition 4.5. Because some of the arguments in the direct treatment can be used again, we do not repeat all details in the proofs of the following propositions.

\iinspr{A.10} Proposition \rm
Denote 
$$F_1=(\iota\ot S)E
\qquad\quad\text{and}\qquad\quad
F_2=(S\ot\iota)E.$$
Then the maps $G_1$ and $G_2$ defined for the pair $(A,\Delta)$ as in Proposition 1.11 satisfy
$$\align 
G_1(a\ot b)&=(a\ot 1)F_1(1\ot b)\\
G_2(a\ot b)&=(a\ot 1)F_2(1\ot b)
\endalign$$
for all $a,b\in A$.

\snl\bf Proof\rm:
We claim that (as in Proposition 4.6)
$$\align 
E_{13}(F_1\ot 1) &= E_{13}(1\ot E) \\
(1\ot F_2)E_{13} &=(E\ot 1)E_{13}.
\endalign$$
To show the first formula, we use that $E(y\ot 1)=E(1\ot S(y))$ for any element in the right leg of $E$ (Lemma A.5) and that $\sigma(S\ot S)E=E$. Then it is a straightforward consequence of the definition of $G_1$ that $G_1(a\ot b)=(a\ot 1)F_1(1\ot b)$ for all $a,b\in A$. 
\snl
Similarly for $G_2$.
\hfill$\square$\einspr

Now, we can continue with showing that the pair $(A^{\text{op}},\Delta)$ satisfies the axioms of a weak multiplier Hopf algebra. Recall that we have already proven that the coproduct $\Delta$ on $A$ also is a full coproduct on the opposite algebra $A^{\text{op}}$ and that the original counit is also the new counit (Proposition A.1). In the next proposition, we see that $E$ is the canonical idempotent for the new pair as well.
\snl
In the remaining of this appendix, in order to avoid confusion, we will use the convention that maps associated with the new pair $(A^{\text{op}},\Delta)$ will be denoted with the same symbols, but with an extra {\it prime}. So e.g.\ the maps $T_1$ and $T_2$ for the new pair will be denoted as $T'_1$ and $T'_2$ respectively. The formulas will always be written down using the original product in $A$.

\iinspr{A.11} Proposition \rm
Consider the maps $T'_1$ and $T'_2$ for the new pair $(A^{\text{op}},\Delta)$. They are given on $A\ot A$ by $T_3$ and $T_4$, defined earlier (after Notation 1.2 in Section 1) by the formulas
$$T_3(a\ot b)=(1\ot b)\Delta(a)
	\qquad\quad\text{and}\qquad\quad
		T_4(a\ot b)=\Delta(b)(a\ot 1).$$
We have 
$$T_3(A\ot A)=(A\ot A)E
	\qquad\quad\text{and}\qquad\quad
		T_4(A\ot A)=E(A\ot A).$$
This shows that the canonical idempotent for the pair $(A^{\text{op}},\Delta)$ is again $E$.

\snl\bf Proof\rm:
From the proof of Proposition A.1, we see that
$$T_3=(\iota\ot S^{-1})R_1(\iota\ot S).$$
We know that the range of $R_1$ is the same as the range of $G_1$. And because 
$$G_1(a\ot b)=(a\ot 1)F_1(1\ot b)$$
for all $a,b$, with $F_1=(\iota\ot S)E$, we get that the range of $T_3$ is precisely $(A\ot A)E$. 
\snl
Similarly from the proof of Proposition A.1, we have
$$T_4=(S^{-1}\ot\iota)R_2(S \ot\iota)$$ 
and as we know that the range of $R_2$ is the same as the range of $G_2$, it follows from the formula
$$G_2(a\ot b)=(a\ot 1)F_2(1\ot b)$$ 
and $F_2=(S\ot 1)E$ that the range of $T_4$ is $E(A\ot A)$.
\snl
Taken into account that the product is reversed for $A^{\text{op}}$, the result is proven.
\hfill$\square$\einspr

We are now ready to finish the proof of the main theorem here.

\iinspr{A.12} Theorem \rm
Let $(A,\Delta)$ be a weak multiplier Hopf algebra. Assume that the antipode $S$ maps $A$ to $A$ and that it is bijective. Then $(A^{\text{op}},\Delta)$ is again a weak multiplier Hopf algebra. In other words, the original pair $(A,\Delta)$ is a {\it regular} weak multiplier Hopf algebra.
\einspr

If we look at the axioms for a weak multiplier Hopf algebra in Definition 1.14, for the new pair $(A^{\text{op}},\Delta)$, we see that all requirements are proven, except for the last property, giving the kernels of the associated maps $T'_1$ and $T'_2$. 
\snl
However, as we have the candidate for the antipode $S'$ for the new pair, namely $S^{-1}$, it turns out to be easier to complete the proof by using the alternative definition, given in Theorem 2.9.

\iinspr{} Proof \rm (of Theorem A.12): 
i) We have already shown that the coproduct $\Delta$ is still a full coproduct on the opposite algebra $A^{\text{op}}$. And of course, $A^{\text{op}}$ is a non-degenerate idempotent algebra. The map $S'$ defined as $S^{-1}$ is here a linear map from $A$ to $A$ and if we define the associated maps $R'_1$ and $R'_2$ given by this antipode $S'$ on the opposite algebra $A^{\text{op}}$, we find as argued before that
$$R'_1=(\iota\ot S^{-1})T_1(\iota\ot S)
 	\qquad\quad\text{and}\qquad\quad
 R'_2=(S^{-1}\ot\iota)T_2(S\ot\iota).$$
 In particular, we see that they have range in $A\ot A$. This is the first requirement in item i) of Theorem 2.9. The formulas (2.5) in this item of Theorem 2.9 for $S'$ on $A^{\text{op}}$ are obtained by simply applying $S^{-1}$ to the same formulas for $S$ on the original algebra $A$.
\snl
ii) Now  we can calculate $T'_1R'_1$ and we find $(\iota\ot S^{-1})R_1T_1(\iota\ot S)$. Because we know that 
$$R_1T_1(a\ot b)=G_1(a\ot b)=(a\ot 1)F_1(1\ot b)$$
and $F_1=(\iota\ot S)E$, we find that 
$$T'_1R'_1(a\ot b)=(a\ot b)E$$
for all $a,b\in A$. Similarly we will have $T'_2R'_2(a\ot b)=E(a\ot b)$ for all $a,b\in A$. This takes care of the first part of item ii) in Theorem  2.9, i.e.\ the formulas (2.6) (taking into account that the product is reversed). And of course, since the idempotent $E'$ for the new pair is the same as for the old pair, also the second part of item ii) in Theorem 2.9, the formulas (2.7), is fulfilled.
\snl
So we can apply Theorem 2.9 and find that $(A^{\text{op}},\Delta)$ is a weak multiplier Hopf algebra.
\hfill$\square$\einspr
 
%
\nl

\bf References \rm
\nl
{[\bf A]} E.\ Abe: {\it Hopf algebras}. \rm Cambridge University Press (1977).
\snl
{[\bf B-N-S]} G.\ B\"ohm, F.\ Nill \& K.\ Szlach\'anyi: {\it Weak Hopf algebras I. Integral theory and C$^*$-structure}. J.\ Algebra 221 (1999), 385-438. 
\snl
{[\bf B-S]} G.\ B\"ohm  \& K.\ Szlach\'anyi: {\it Weak Hopf algebras II. Representation theory, dimensions and the Markov trace}. J.\ Algebra 233 (2000), 156-212. 
\snl
{[\bf Br]} R.\ Brown: {\it From groups to groupoids: A brief survey}. Bull. London Math. Soc. 19 (1987), 113–134.
\snl
{[\bf D-VD-Z]} B.\ Drabant, A.\ Van Daele \& Y.\ Zhang: {\it Actions of multiplier Hopf algebras}. Commun.\ in alg.\ 27 (1999), 4117-4172.
\snl
{[\bf G]} K.R.\ Goodearl: {\it Von Neumann Regular Rings}. Pitman, London, 1979.
\snl
{[\bf H]} P.\ J.\ Higgins: {\it Notes on categories and groupoids}. Van Nostrand Reinhold, London (1971).
\snl
{[\bf N]} D.\ Nikshych: {\it On the structure of weak Hopf algebras}. Adv.\ Math.\ 170 (2002), 257-286.
\snl
{[\bf N-V1]} D.\ Nikshych \& L.\ Vainerman: {\it Algebraic versions of a finite dimensional quantum groupoid}. Lecture Notes in Pure and Applied Mathematics 209 (2000), 189-221. 
\snl
{[\bf N-V2]} D.\ Nikshych \& L.\ Vainerman: {\it Finite quantum groupiods and their applications}. In {\it New Directions in Hopf algebras}. MSRI Publications, Vol.\ 43 (2002), 211-262.
\snl
{[\bf P]} A.\ Paterson: {\it Groupoids, inverse semi-groups and their operator algebras}. Birkhauser, Boston (1999).
\snl
{[\bf R]} J.\ Renault: {\it A groupoid approach to C$^*$-algebras}. Lecture Notes in Mathematics 793, Springer Verlag.
\snl
{[\bf S]} M.\ Sweedler: {\it Hopf algebras}. Benjamin, New-York (1969).
\snl
{[\bf V]} L.\ Vainerman (editor): {\it Locally compact quantum groups and groupoids}. IRMA Lectures in Mathematics and Theoretical Physics 2, Proceedings of a meeting in Strasbourg, de Gruyter (2002).
\snl
{[\bf VD1]} A.\ Van Daele: {\it Multiplier Hopf algebras}. Trans. Am. Math. Soc.  342(2) (1994), 917-932.
\snl
{[\bf VD2]} A.\ Van Daele: {\it An algebraic framework for group duality}. Adv. in Math. 140 (1998), 323-366.
\snl
{[\bf VD3]} A.\ Van Daele: {\it Tools for working with multiplier Hopf algebras}. Preprint University of Leuven. Arxiv math.QA/0806.2089. ASJE (The Arabian Journal for Science and Engineering) C - Theme-Issue 33 (2008), 505--528.  
\snl
{[\bf VD4]} A.\ Van Daele: {\it Separability and multiplier algebras}. Preprint University of Leuven (2012) (in preparation).
\snl
{[\bf VD-Ve]} A.\ Van Daele \& J.\ Vercruysse: {\it Local units and multiplier algebras}. Preprint University of Leuven and University of Brussels (2012) (in preparation).
\snl
{[\bf VD-W1]} A.\ Van Daele \& S.\ Wang: {\it The Larson-Sweedler theorem for multiplier Hopf algebras}. J.\ of Alg.\ 296 (2006), 75-95.
\snl
{[\bf VD-W2]} A.\ Van Daele \& S.\ Wang: {\it Multiplier Unifying Hopf Algebras}. Preprint University of Leuven and Southeast University of Nanjing (2008) (unpublished).
\snl 
{[\bf VD-W3]} A.\ Van Daele \& S.\ Wang: {\it Weak multiplier Hopf algebras. Preliminaries, motivation and basic examples}. Preprint University of Leuven and Southeast University of Nanjing (2012). Arxiv:1210.3954v1 [math.RA]. To appear in the proceedings of  the conference 'Operator Algebras and Quantum Groups (Warsaw, September 2011), series 'Banach Center Publications'.
\snl 
{[\bf VD-W4]} A.\ Van Daele \& S.\ Wang: {\it Weak multiplier Hopf algebras II. The source and target algebras}. Preprint University of Leuven and Southeast University of Nanjing (2012) (in preparation).
\snl
{[\bf VD-W5]} A.\ Van Daele \& S.\ Wang: {\it Weak multiplier Hopf algebras III. Integrals and Duality}. Preprint University of Leuven and Southeast University of Nanjing (in preparation).
\snl
{[\bf Ve]} J.\ Vercruysse: {\it Local units versus local projectivity dualisations: Corings with local structure maps}. Commun. in Alg. 34 (2006) 2079-2103.

\end